\documentclass[11pt]{article}

\usepackage{tikz}
\usepackage{subfigure}
\usepackage{hyperref}
\usepackage[text={450pt,600pt}]{geometry}
\usepackage{titlesec}
\usepackage{amssymb}
\usepackage{xcolor}
\usepackage{amsmath,amssymb,amscd,mathrsfs}
\usepackage{color}
\usepackage{amssymb}
\usepackage{amsmath}
\usepackage{amscd}
\usepackage{ntheorem}
\usepackage{xcolor}

\DeclareMathOperator{\Prog}{Prog}

 \newcommand{\sgn}{\operatorname{sgn}}
 
\usepackage{cleveref}

\hypersetup{
  colorlinks   = true,    
  urlcolor     = blue,    
  linkcolor    =  blue,  
  citecolor    = blue      
}
\newcommand*{\qed}{\par\vspace{5mm} \hfill\ensuremath{\Box}}
\usepackage[latin1]{inputenc}   
\usepackage[T1]{fontenc}   
\usepackage{mathtools}
\DeclarePairedDelimiterX{\norm}[1]{\lVert}{\rVert}{#1}

\begin{document}

\renewcommand\thesubsection{\arabic {subsection}}
\newtheorem{thm}{Theorem} [section]
  \newtheorem{Conjecture}[thm]{Conjecture}
   \newtheorem{Lemma}[thm]{Lemma} 
\newtheorem{prop}[thm]{Proposition}
 \newtheorem{Corollary}[thm]{Corollary}
     \newtheorem{ex}{Example}
    \theorembodyfont{\normalfont}
  \newtheorem{rem}{Remark} 
     \newtheorem{as}{Assumption}[section]
\renewcommand{\theas}{{{(A.\arabic{as}})}}
\newtheorem{defi}{Definition}[section]
  \renewcommand{\thedefi}{\arabic{defi}.}
  \renewcommand{\therem}{\arabic{rem}.}
   \theorembodyfont{\normalfont}\theoremstyle{nonumberplain}  
     \newtheorem{Proof}{Proof.}
  \newtheorem{remm}{Remark.} 

  \renewcommand{\thethm}{\arabic{thm}}%

\newenvironment{thmbis}[1]
  {\renewcommand{\thethm}{\ref{#1}$'$}%
   \addtocounter{thm}{-1}%
   \begin{thm}}
  {\end{thm}}

\newenvironment{asp}[1]
 {\addtocounter{as}{-1}%
  \renewcommand{\theas}{(A.\arabic{as}$'$)}%
   \begin{as}}
  {\end{as}}

\newenvironment{aspp}[1]
 {\addtocounter{as}{-1}%
  \renewcommand{\theas}{(A.\arabic{as}$''$)}%
   \begin{as}}
  {\end{as}}

\makeatletter
\newcommand{\neutralize}[1]{\expandafter\let\csname c@#1\endcsname\count@}
\makeatother
\title{\bf Lipschitz-quadratic Regularization for Quadratic Semimartingale BSDEs}
\author{Hanlin Yang\thanks{{Department of Banking and Finance, Universit\"at Z\"
urich, Plattenstrasse 22, CH-8006 Z\"urich, Switzerland, e-mail: hanlin.yang@uzh.ch. This paper is part of the author's master thesis.}}} 
\date{\today}
\maketitle
\abstract{
We refine the solvability of quadratic semimartingale BSDEs by employing a  Lipschitz-quadratic regularization procedure. In the first step, we prove an existence and uniqueness result for a class of Lipschitz-quadratic BSDEs. 
A corresponding  stability theorem and a Lipschitz-quadratic regularization
are developed to solve quadratic BSDEs. The advantage of our approach is that much weaker conditions  ensure the existence and uniqueness results.
\subsection{Preliminaries}
\label{section20}
In this paper, we are concerned with the solvability of $\mathbb{R}$-valued backward stochastic differential equations (BSDEs) driven by continuous local martingales which take the form
\begin{align}
\label{eqn1}
Y_t = \xi + \int_t^T\big( \mathbf{1}^\top
d\langle M\rangle_s F(s, Y_s, Z_s)
+ g_s d\langle N\rangle_s\big) -\int_t^T \big(Z_s dM_s + dN_s\big),
\end{align}
where $M$ and $N$ are strongly orthogonal continuous local martingales. We are particularly interested in the above equations with quadratic growth, i.e., the generator $F$ is quadratic in $Z$ and $g$ is not identical to $0$.

BSDEs of this type have been intensively applied to mathematical finance and stochastic control; see Mania and Schweizer \cite{MS2005}, or Hu et al \cite{HIM2005} in Brownian setting. In its theoretical aspect,  Karoui and Huang \cite{K1994} obtains the solvability with Lipschitz-continuous generators. 
  Later, Tevzadze \cite{T2008} studies the existence and uniqueness of a bounded solution, by assuming quadratic growth and local Lipschitz-continuity.
 Morlais \cite{M2009} extends the  stability-type argument  in Kobylanski \cite{K2000} to quadratic BSDEs driven by continuous local martingales. Based on this work, Mocha and Westray 
 \cite{MW2012} proves existence and uniqueness results with convex generators and exponential moments integrability.

A close inspection of this line of study, however,  reveals that their assumptions are quire demanding. 
 For example, the stability-type argument in Morlais \cite{M2009}  can be used only if the BSDE is not quadratic in $N$, i.e., $g_\cdot = 0$. When $g$ is a  constant process, an exponential transform can be used to kill the quadratic term $g\cdot \langle N\rangle$. But one has to sacrifice the flexibility of the generators, especially for unbounded solutions.
 For this point, the interested readers shall refer to \cite{K1994}, \cite{M2009}, \cite{MW2012}.
When $g$ is a bounded process, some results are obtained by Tevzdaze \cite{T2008}, but rather restrictive. For example, existence results are obtained only for  particular quadratic generators, and equations with Lipschitz-continuous  generators are not studied.

Having understood these literature and their drawbacks, we develop a Lipschitz-quadratic regularization technique to answer the question of existence and uniqueness under more flexible assumptions. In the first step,  we study BSDEs with  Lipschitz-continuous generators and quadratic growth in  $N$, by adapting the fixed point arguments in Tevzadze \cite{T2008}.  These equations, due to this particular structure, are called Lipschitz-quadratic.
Viewing this result as a basic building block, we  then derive a corresponding monotone stability result to faciliate our study of more general quadratic BSDEs. The  regularization therein is called Lipschitz-quadratic, as contrary to the Lipschitz regularization in \cite{M2009}, \cite{MW2012}. It turns out that all the results, including existence, uniqueness and stability results of bounded and unbounded solutions 
can be obtained with  weaker conditions.

This paper is organized as follows. In Section \ref{section21}, we prove an existence and uniqueness result for Lipschitz-quadratic BSDEs. Based on this result, we establish a monotone stability theorem in Section \ref{section22}. As a byproduct, the existence of a bounded solution is immediate.
 In Section \ref{section23}, we  study existence, uniqueness  and stability results for unbounded solutions, using a localization procedure.  Finally, Section \ref{section24} reviews the change of measure result studied in Mocha and Westray 
 \cite{MW2012}.

Let us close this section by introducing all required notations.
  We fix the time horizon $0<T<+\infty$, and work on a filtered probability space $(\Omega, \mathcal{F}, (\mathcal{F}_t)_{t\in [0, T]}, \mathbb{P})$ satisfying the usual conditions of right-continuity and $\mathbb{P}$-completeness. 
 $\mathcal{F}_0$ is the $\mathbb{P}$-completion of the trivial $\sigma$-algebra.
Any measurability  will refer to the filtration $(\mathcal{F}_t)_{t\in [0, T]}$.  In particular,  $\Prog$ denotes the progressive $\sigma$-algebra on $\Omega \times [0, T]$.
We assume the filtration  is \emph{continuous}, in the sense that all  local martingales have $\mathbb{P}$-a.s. continuous sample paths.  
 $M = (M^1,..., M^d)^\top$  stands for a  fixed $d$-dimensional   continuous local martingale.  By \emph{continuous semimartingale setting} we mean:  $M$ doesn't have to be a Brownian motion;  the filtration is not necessarily generated by $M$ which is usually seen as the main source of randomness. Hence in various concrete situations there may be a continuous  local martingale strongly orthogonal to $M$, which we denote, as in (\ref{eqn1}), by $N$.

Here we clarify all notions in (\ref{eqn1}). We set $\mathbf{1}:=(1, ..., 1)^\top$. $\xi$ is an $\mathbb{R}$-valued $\mathcal{F}_T$-measurable random variable, $F: \Omega \times [0, T]
\times \mathbb{R} \times \mathbb{R}^d \rightarrow \mathbb{R}^d$ is a $\Prog\otimes \mathcal{B}(\mathbb{R})\otimes
\mathcal{B}(\mathbb{R}^d)$-measurable random function and $g$ is an $\mathbb{R}$-valued $\Prog$-measurable bounded process.  $\int_0^\cdot (Z_sdM_s + dN_s)$, sometimes denoted by $Z\cdot M + N$,
refers to the vector stochastic integral; see Shiryaev and Cherny \cite{SC2002}.
The equations defined in this way encode the matrix-valued process $\langle M\rangle$  which is not amenable to analysis. Therefore we rewrite the BSDEs by factorizing
$\langle M \rangle$. This procedure separates the matrix property from its nature as a measure. It can also be regarded as a reduction of dimensionality.

There are many ways to factorize $\langle M \rangle$; see, e.g., Section III. 4a, Jacod and  Shiryaev \cite{JS1987}. We can and  choose
$A: = \arctan\big( \sum_{i=1}^d {\langle M^i\rangle}\big)$.   By Kunita-Watanabe inequality, we deduce the absolute continuity of  $\langle M^i, M^j \rangle$  
with respect to $A$. Note that such choice makes $A$ continuous, increasing and bounded.
Moreover, by  Radon-Nikod\'{y}m theorem and Cholesky decomposition, there exists a matrix-valued $\Prog$-measurable process $\lambda$  
such that $\langle M\rangle  = (\lambda^\top \lambda) \cdot A.$ 
As will be seen later, our results
don't rely on the specific choice of $A$ but only on its boundedness. In particular, if $M$ is a $d$-dimensional Brownian motion, we may choose $A_t =t$ and $\lambda$ to be the  identity matrix.

 The second advantage of  factorizing $\langle M\rangle$  is that
\[
\mathbf{1}^\top d\langle M\rangle_s F(s, Y_s, Z_s) =
\mathbf{1}^\top \lambda_s^\top \lambda_s F(s, Y_s, Z_s)dA_s,
\]
where $f(t, y, z) := \mathbf{1}^\top \lambda_s^\top \lambda_s F(s, y, z)$  is $\mathbb{R}$-valued. Such reduction of dimensionality makes it easier to formulate the difference of two equations as frequently appears in comparison theorem and uniqueness.
Hence, we may reformulate the BSDEs as follows.

{\bf BSDEs: Definition and Solutions.}
Let $A$ be an $\mathbb{R}$-valued continuous  nondecreasing {bounded} adapted process such that 
$\langle M \rangle = (\lambda^\top \lambda) \cdot A$ for some matrix-valued $\Prog$-measurable process $\lambda$, $f: \Omega \times [0, T]
\times \mathbb{R} \times \mathbb{R}^d \rightarrow \mathbb{R}$
  a $\Prog\otimes \mathcal{B}(\mathbb{R})\otimes
\mathcal{B}(\mathbb{R}^d)$-measurable random function,
$g$  an $\mathbb{R}$-valued $\Prog$-measurable bounded process and $\xi$  an $\mathbb{R}$-valued $\mathcal{F}_T$-measurable random variable.  The semimartingale BSDEs are written as
\begin{align}
Y_t= \xi + &\int_t^T \big(f(s, Y_s, Z_s)dA_s + g_s d\langle N\rangle_s\big) -
\int_t^T \big(Z_s dM_s + dN_s\big).  \label{sbsde}
\end{align}
 We call a process $(Y, Z, N)$   or $(Y, Z\cdot M + N)$ a \emph{solution} of (\ref{sbsde}), if $Y$ is an $\mathbb{R}$-valued continuous adapted process, $Z$ is an $\mathbb{R}^d$-valued $\Prog$-measurable process and $N$  is 
an $\mathbb{R}$-valued 
  continuous local martingale strongly orthogonal to $M$, such that $\mathbb{P}$-a.s. $\int_0^T Z_s^\top d\langle M\rangle_s Z_s <+\infty$ and 
$\int_0^T |f(s, Y_s, Z_s)|dA_s <+\infty$, and  (\ref{sbsde}) holds   
 $\mathbb{P}$-a.s. for all $t\in [0, T]$,

Note that the factorization of $\langle M\rangle$ gives $
 \int_0^\cdot Z_s^\top d\langle M\rangle_s Z_s  = \int_0^\cdot |\lambda_s Z_s|^2 dA_s
$ $\mathbb{P}$-a.s.
Hence we don't  distinguish these two integrals  in all situations. 
$\int_0^T Z_s^\top d\langle M\rangle_s Z_s <+\infty\ \mathbb{P}$-a.s. ensures  that $Z$ is integrable with respect to $M$ in the sense of vector stochastic integration.
As a result, $Z\cdot M$ is a continuous local martingale. 
$M$ and $N$ being continuous and strongly orthogonal implies that
$\langle M^i, N\rangle_\cdot = 0$ for  $i = 1, ..., d$.
We call $f$ the \emph{generator}, $\xi$ the \emph{terminal value} and   $(\xi, \int_0^T|f(s, 0, 0)|dA_s)$ the \emph{data}. 
In our study, the integrability property of the data determines the estimates for a solution. The conditions imposed on the generator are called the \emph{structure conditions}. For notational convenience, we sometimes write $(f, g, \xi)$ instead of (\ref{sbsde}) to denote the above BSDE.  Finally, 
(\ref{sbsde})
 is called \emph{quadratic} if $f$ has at most quadratic growth in $z$ or $g$ is not indistinguishable from $0$.
 
To finalize, we introduce the rest notations  which will be used throughout this paper.
 $\ll$ stands for the strong order of nondecreasing processes, stating that the difference is nondecreasing.
For any random variable or process $Y$, we say $Y$ has some property if this is true except on a $\mathbb{P}$-null subset of $\Omega$. Hence we omit ``$\mathbb{P}$-a.s'' in situations without ambiguity. Define $\sgn(x)= \mathbb{I}_{\{x\neq 0 \}}\frac{x}{|x|}$. For any random variable $X$, define $\norm{X}_\infty$ to be its essential supremum.
For any c\`adl\`ag adapted process $Y$, set $Y_{s, t}: = Y_t -Y_s$ and $Y^* : = \sup_{t\in [0, T]} |Y_t|$.
For any $\Prog$-measurable process $H$, set $|H|_{s,t}:= \int_s^t H_u dA_u$ and $|H|_t : = |H|_{0, t}$. 
 $\mathcal{T}$ stands for the set of all stopping times valued in $[0, T]$ and $\mathcal{S}$ denotes the space of continuous adapted processes. 
For later use we specify the following spaces under $\mathbb{P}$.
\begin{itemize}
\item $\mathcal{S}^\infty$: the space of  bounded processes  $Y\in\mathcal{S}$ with $\norm{Y}:=\norm{Y^*}_\infty$; $\mathcal{S}^\infty$ is a  Banach space;

	\item $\mathcal{M}$: the set of  continuous local martingales starting from $0$; for any $\mathbb{R}^d$-valued $\Prog$-measurable process $Z$ with $\int_0^T Z_s^\top d\langle M\rangle_s Z_s < +\infty$,  $Z\cdot M \in \mathcal{M}$; 
\item $\mathcal{M}^p(p\geq 1)$: the set 
of $\widetilde{M} \in \mathcal{M}$ with 
\[
\norm{\widetilde{M}}_{\mathcal{M}^p}:=\big(\mathbb{E}\big[ \langle \widetilde{M}\rangle_T^{\frac{p}{2}}\big] \big)^{\frac{1}{p}} < +\infty;
\]
 in particular, $\mathcal{M}^2$ is a Hilbert space;
\item $\mathcal{M}^{BMO}$: the 
 set  of BMO martingales $\widetilde{M}\in \mathcal{M}$ \ with 
\[ \norm{\widetilde{M}}_{BMO} :=\sup_{\tau \in \mathcal{T}}
\big\|\mathbb{E}\big[
\langle \widetilde{M}\rangle_{\tau, T}
\big|\mathcal{F}_\tau    \big]^{\frac{1}{2}}\big\|_\infty;
\]
$\mathcal{M}^{BMO}$ is a Banach space.
\end{itemize}

 $\mathcal{M}^2$ being a Hilbert space is crucial to proving convergence of the  martingale parts in the monotone stability result of quadratic BSDEs; see, e.g.,  Kobylanski \cite{K2000},  Briand and  Hu \cite{BH2008},  Morlais \cite{M2009} or Section \ref{section22}. 
 Other spaces  are also Banach  under suitable norms; we will not present these facts in more detail since they are not involved in our study.
  
Finally, for any local martingale $\widetilde{M}$, we call 
$\{\sigma_n \}_{n\in\mathbb{N}^+}\subset \mathcal{T}$ a \emph{localizing sequence} if $\sigma_n$ increases stationarily to $T$ as $n$  goes to $+\infty$ and $\widetilde{M}_{\cdot\wedge \sigma_n}$ is a martingale for any $n\in\mathbb{N}^+$.

\subsection{Bounded Solutions of Lipschitz-quadratic BSDEs}
\label{section21}
This section takes one step in   solving   quadratic BSDEs and consists  in the study of equations with Lipschitz-continuous generators. In contrast to  El Karoui and  Huang \cite{EK1997}, we allow the presence of  $g\cdot \langle N\rangle$.
We point out that similar results
for linear-quadratic generators have been studied by
Tevzadze \cite{T2008}, but the case of Lipschitz-continuity is not available in that work.   Due to its importance for regularizations of quadratic BSDEs, we study existence and uniqueness   results for equations of this particular type  in the first step.  To this end, we assume
\begin{as}\label{as1} 
There exist $\beta, \gamma \geq 0$
such that 
${\norm{\xi}}_\infty + \big\|{\big| |f(\cdot, 0, 0) | \big|_T}\big\|_\infty <+\infty$ and 
$f$ is Lipschitz-continuous in $(y, z)$, i.e., $\mathbb{P}$-a.s.  for any $t\in[0, T]$, $y, y^\prime \in\mathbb{R}$, $z, z^\prime \in \mathbb{R}^d$,
\[
|f(t, y, z ) - f(t, y^\prime, z^\prime)| \leq \beta |y-y^\prime| + \gamma |\lambda_t (z-z^\prime)|.
\]
\end{as}

Due to the presence of $g\cdot\langle N\rangle$,  we call the BSDE $(f, g, \xi)$ satisfying \ref{as1}   \emph{Lipschitz-quadratic}. Given \ref{as1}, 
we are about to construct a solution in the  space $\mathscr{B}: = \mathcal{S}^\infty\times \mathcal{M}^{BMO}$ equipped with  the norm 
\[
\norm {(Y, Z\cdot M + N)} : = \big(\norm{Y}^2 + \norm{Z\cdot M + N }^2_{BMO}\big)^\frac{1}{2},
\] 
for $(Y, Z\cdot M  + N)\in \mathcal{S}^\infty\times\mathcal{M}^{BMO}$. Clearly $(\mathscr{B}, \norm{\cdot})$ is Banach. 
As a preliminary result,  we claim that   the existence result holds given sufficiently small data.
\begin{thm}[Existence (i)]\label{existence1} If $(f, g, \xi)$ satisfies  {\rm\ref{as1}} with 
\begin{align}
\norm{\xi}^2_{\infty}+8 \big\|{\big||f(\cdot, 0, 0)|\big|_T}\big\|^2_{\infty} \leq \frac{1}{64}\exp \Big(       -\norm{A} 
\big(8\beta^2 \norm{{A}} + 8\gamma^2\big)
    \Big) \label{p1}
\end{align}
 and $\mathbb{P}$-a.s.  $|g_\cdot|\leq \tilde{g}:= \frac{1}{8}$, then
 there exists a solution  in $(\mathscr{B}, \norm{\cdot})$.
\end{thm}	
\begin{Proof}  To overcome the difficulty arising from  the Lipschitz-continuity, we use Banach fixed point theorem under an equivalent norm.
Set $\rho \geq 0$ to be determined later. 
For any  $X \in\mathbb{L}^\infty, Y\in\mathcal{S}^\infty$ and $\widetilde{M} \in\mathcal{M}^{BMO}$, set
$\norm{X}_{ \infty, \rho} : = 
\norm{e^{\frac{\rho}{2} A_T}X }_\infty$,
 $\norm{Y}_\rho : = \norm{e^{\frac{\rho}{2}A} Y } $ and $\norm{\widetilde{M}}_{BMO, \rho} : = \norm{e^{\frac{\rho}{2} A}\cdot \widetilde{M}}_{BMO}$;   for $(Y, Z\cdot M +N)\in \mathscr{B}$, set
\[
\norm{(Y, Z\cdot M +N)}_\rho :  =  \big(\norm{Y}_\rho^2 + \norm{Z\cdot M + N}_{BMO, \rho}^2\big)^{\frac{1}{2}}.
\]
Since $A$ is bounded, $\norm{\cdot}_\rho$  is equivalent to the original norm for each space. Hence
$(\mathscr{B}, \norm{\cdot}_\rho)$ is also a Banach space.
For any $R\geq 0$,  define
\[
\mathbf{B}_R : = \big\{ (Y, Z\cdot M + N) \in \mathscr{B}: \norm{(Y, Z\cdot M +N)}_\rho \leq R   \big\}.
\]
We show by  Banach fixed point theorem that there exists a unique solution in $\mathbf{B}_{R}$ with $R=\frac{1}{2}.$
To this end, we define $\mathbf{F}: (\mathbf{B}_R, \norm{\cdot}_\rho)\rightarrow (\mathscr{B},  \norm{\cdot}_\rho)
$ such that for any $(y, z\cdot M + n) \in\mathbf{B}_R$,  $(Y, Z\cdot M + N) : =\mathbf{F}((y, z\cdot M + n))$ solves
\[
Y_t = \xi + \int_t^T \big( f(s, y_s, z_s) dA_s + g_s d\langle n \rangle_s \big) -\int_t^T \big(Z_s dM_s + dN_s \big).
\]
Indeed, such $(Y, Z, N)$  uniquely exists due to martingale representation theorem. Moreover, by standard estimates, $(Y, Z\cdot M + N) \in (\mathscr{B}, \norm{\cdot}_\rho)$.

(i). We show
$
\mathbf{F}(\mathbf{B}_R) \subset \mathbf{B}_R.
$
For any $\tau \in \mathcal{T}$,  It\^o's formula applied to $e^{\rho A_\cdot}Y^2_\cdot$ yields
\begin{align}
&e^{\rho A_\tau}|Y_\tau|^2 + \rho \mathbb{E}\Big[ \int_\tau^T e^{\rho A_s}Y_s^2 dA_s \Big|\mathcal{F}_\tau \Big]
+ \mathbb{E}\Big[\int_\tau^T e^{\rho A_s} \Big(Z_s^\top d\langle M\rangle_s  Z_s + d\langle N\rangle_s\Big)\Big|\mathcal{F}_\tau \Big] 
\nonumber\\ 
&\leq
\norm{\xi}^2_{\infty, \rho}  + 2\mathbb{E}\Big[\int_\tau^T e^{\rho A_s}|Y_s ||f(s, y_s, z_s)|dA_s \Big| \mathcal{F}_\tau\Big] + 2\mathbb{E}\Big[ \int_\tau^T e^{\rho A_s}|Y_s||g_s| d\langle n\rangle_s\Big| \mathcal{F}_\tau\Big].
\label{lq1}
\end{align}
By  \ref{as1},
\begin{align*}
|Y_s||f(s, y_s, z_s)| \leq |Y_s||f(s, 0, 0)| + \beta |Y_s| |y_s| + \gamma |Y_s| |\lambda_s z_s|.
\end{align*}
We plug this inequality into (\ref{lq1})
and estimate each term on the right-hand side.  Using $2ab \leq \frac{1}{8}a^2 + 8b^2$  gives  
\begin{align*}
2\mathbb{E}\Big[  \int_\tau^T  e^{\rho A_s} |Y_s||f(s, 0, 0)| dA_s     \Big| \mathcal{F}_\tau\Big]
&\leq 
\frac{1}{8}\norm{Y}^2_\rho +  8 \mathbb{E}\Big[ \int_\tau^T e^{\frac{\rho}{2}A_s} |f(s, 0, 0)|dA_s  \Big|\mathcal{F}_\tau\Big]^{2} \\
&\leq \frac{1}{8}\norm{Y}^2_\rho + 8 \big\|{\big||f(\cdot, 0, 0)|\big|_T}\big\|^2_{\infty, \rho},
\end{align*}
\begin{align*}2
\beta \mathbb{E}\Big[ \int_\tau^T e^{\rho A_s}|Y_s||y_s| dA_s \Big|\mathcal{F}_\tau\Big]
&\leq
\frac{1}{8}\norm{y}^2_\rho + 8\beta^2 \mathbb{E}\Big[\int_\tau^T e^{\frac{\rho}{2}A_s} |Y_s| dA_s    \Big| \mathcal{F}_\tau \Big]^{2}\\
& \leq \frac{1}{8}\norm{y}^2_\rho + 8\beta^2\norm{A} \mathbb{E}\Big[\int_\tau^T e^{
\rho A_s} |Y_s|^2 dA_s    \Big| \mathcal{F}_\tau \Big], 
\end{align*}
\begin{align*}
2\gamma \mathbb{E}\Big[ \int_\tau^T e^{\rho A_s}|Y_s||\lambda_s z_s| dA_s \Big|\mathcal{F}_\tau\Big]
& \leq \frac{1}{8}\norm{z\cdot M}^2_{BMO, \rho} + 8\gamma^2 \mathbb{E}\Big[\int_\tau^T e^{\rho A_s} |Y_s|^2 dA_s    \Big| \mathcal{F}_\tau \Big],
\end{align*}
\begin{align*}
2 \mathbb{E}\Big[\int_\tau^T e^{\rho A_s }|Y_s| |g_s|\langle N \rangle _s \Big| \mathcal{F}_\tau \Big]
&\leq \frac{1}{8}\norm{Y}^2_\rho + 8 \tilde{g}^2\mathbb{E}\Big[ \int_\tau^T e^{\frac{\rho}{2}A_s} d\langle N\rangle_s \Big| \mathcal{F}_\tau        \Big]^2 \\
&\leq \frac{1}{8}\norm{Y}^2_\rho + 8 \tilde{g}^2 \norm{n}^4_{BMO, \rho}.
\end{align*}
Set $\rho :=  8\beta^2 \norm{A} + 8\gamma^2$  so as to eliminate 
$\mathbb{E}\big[ \int_\tau^T e^{\rho A_s}Y_s^2 dA_s \big|\mathcal{F}_\tau \big]$ on both sides. Hence  (\ref{lq1}) gives
\begin{align}
e^{\rho A_\tau}|Y_\tau|^2 
&+ \mathbb{E}\Big[\int_\tau^T e^{\rho A_s} \Big(Z_s^\top d\langle M\rangle_s  Z_s + d\langle N\rangle_s\Big)\Big|\mathcal{F}_\tau \Big]  \nonumber \\
& \leq 
\norm{\xi}^2_{\infty, \rho} +8  \big\| \big| |f(\cdot, 0, 0)|\big|_T\big\|^2_{\infty, \rho} + \frac{1}{4} \norm{Y}^2_\rho  \nonumber\\
&+ \frac{1}{8} \big(\norm{y}^2_\rho + \norm{z\cdot M}_{BMO, \rho}^2 \big)
+ 8 \tilde{g}^2 \norm{n}^4_{BMO, \rho}. \label{lq11}
\end{align}
Taking essential supremum and  supremum over all $\tau \in \mathcal{T}$, and using the inequality
\begin{align*}
\frac{1}{2}\norm {(Y, Z\cdot M +N)}_\rho^2 &\leq   \norm{Y}^2_\rho \vee \norm{Z\cdot M + N}^2_{BMO, \rho}\\
 &\leq  \sup_{\tau \in \mathcal{T}}\Big\| e^{\rho A_\tau}|Y_\tau|^2 +  \mathbb{E}\Big[ \int_\tau^T e^{\rho A_s} \Big(Z_s^\top d\langle M\rangle_s Z_s + d\langle N \rangle_s \Big)\Big| \mathcal{F}_\tau\Big]\Big\|_\infty, 
\end{align*}
we deduce by transferring $\frac{1}{4} \norm{Y}^2_\rho$ to the left-hand side of (\ref{lq11}) that 
\begin{align*}
\norm{(Y, Z\cdot M + N)}^2_\rho &\leq  4\norm{\xi}_{\infty, \rho}^2 + 32  \big\|{\big||f(\cdot, 0, 0)|\big|_T}\big\|^2_{\infty, \rho}+ \frac{1}{2} \big(\norm{y}_\rho^2 + \norm{z\cdot M}^2_{BMO, \rho} \big) +
32  \tilde{g}^2 \norm{n}^4_{BMO, \rho} \\
& \leq 4\norm{\xi}^2_{\infty, \rho} + 32  \big\|{\big||f(\cdot, 0, 0)|\big|_T}\big\|^2_{\infty, \rho}  + \frac{1}{2}R^2  + 32 \tilde{g}^2 R^4. 
\end{align*}
Thanks to (\ref{p1}), $\tilde{g}=\frac{1}{8}$ and $R=\frac{1}{2}$,  we verify from the above estimate that 
\[
\norm{(Y, Z, N)}_\rho \leq R.
\]

(ii).
We  prove $\mathbf{F}: (\mathbf{B}_R, \norm{\cdot}_\rho) \rightarrow (\mathbf{B}_R, \norm{\cdot}_\rho)$ is a contraction mapping. 
By (i),  for $i=1, 2$ and  any $(y^i, z^i\cdot M + n^i)\in \mathbf{B}_R$,   we have  $(Y^i, Z^i\cdot M + N^i) : = \mathbf{F}((y^i, z^i\cdot M + n^i)) \in \mathbf{B}_R$. 
For notational convenience we set 
 $\delta y: = y^1 - y^2$ and  $\delta z, \delta n, \delta \langle n\rangle, \delta Y, \delta Z, \delta N, \delta \langle N\rangle$, etc.  analogously.
By  the deductions  in (i) with  minor modifications, we obtain  
\begin{align}
\frac{1}{2}\norm{(\delta Y, \delta Z\cdot M + \delta N)}^2_\rho &\leq   \frac{1}{8} \big(\norm{\delta y}^2_\rho + \norm{\delta z\cdot M}^2_{BMO, \rho} \big)   + \frac{1}{4}\norm{\delta Y}^2_\rho \nonumber\\
&+ 
4\tilde{g}^2  \sup_{\tau \in \mathcal{T}} \Big \| \mathbb{E} \Big[\int_\tau^T e^{\frac{\rho}{2}A_s} d |\delta \langle n \rangle _s|              \Big| \mathcal{F}_\tau \Big]^2\Big\|_\infty.\label{lq111}
\end{align}
 Kunita-Watanabe inequality and Cauchy-Schwartz inequality  used to the last term gives
\begin{align*}
\mathbb{E} \Big[\int_\tau^T e^{\frac{\rho}{2}A_s} d |\delta \langle n \rangle _s|              \Big| \mathcal{F}_\tau \Big]^2 &\leq
\mathbb{E}\Big[  \int_\tau^T  e^{\frac{\rho}{2}A_s}d\langle \delta n \rangle_s  \Big| \mathcal{F}_\tau    \Big]
 \mathbb{E}\Big[  \int_\tau^T  e^{\frac{\rho}{2}A_s}d\langle n^1 + n^2 \rangle_s  \Big| \mathcal{F}_\tau    \Big]
\\
& \leq  \norm{\delta n}^2_{BMO, \rho}  \cdot  2 \big(\norm{n^1}_{BMO, \rho}^2 + \norm{n^2}_{BMO, \rho}^2\big) \\
&\leq \norm{\delta n}^2_{BMO, \rho} \cdot 4R^2,
\end{align*}
where the last inequality is due to $\norm{(y^i, z^i\cdot M + n^i)}_\rho \leq R$, $i = 1, 2.$
Hence (\ref{lq111}) gives
\begin{align*}
\norm{(\delta Y, \delta Z\cdot M + \delta N)}^2_\rho &\leq \frac{1}{2} \big(\norm{\delta y}^2_\rho + \norm{\delta z\cdot M}^2_{BMO, \rho} \big) +  64 \tilde{g}^2 R^2 \norm{\delta n}^2_{BMO, \rho} \\
&\leq \Big(\frac{1}{2} + 64 \tilde{g}^2 R^2\Big) \norm{(\delta y, \delta z\cdot M + \delta n)}^2_\rho \\
&\leq \frac{3}{4}\norm{(\delta y, \delta z \cdot M + \delta n)}^2_\rho,
\end{align*}
i.e.,  $\mathbf{F}: (\mathbf{B}_R, \norm{\cdot}_\rho) \rightarrow (\mathbf{B}_R, \norm{\cdot}_\rho)$ is a contraction mapping. The existence of a solution in $\mathbf{B}_R$ thus follows immediately from Banach fixed point theorem. Finally, since $\norm{\cdot}$ is equivalent to $\norm{\cdot}_\rho$ for $\mathscr{B}$,  the solution also belongs to $(\mathscr{B}, \norm{\cdot})$.
\qed
\end{Proof}

From now on we denote $(\mathscr{B},\norm{\cdot})$ by $\mathscr{B}$ when there is no ambiguity.
 In the spirit of  Tevzadze \cite{T2008},  we  extend this existence result so as to allow any bounded data. 
 To this end, for any $\mathbb{Q}$ equivalent to $\mathbb{P}$ we define 
$\mathcal{S}^\infty (\mathbb{Q})$ analogously to $\mathcal{S}^\infty$ but under $\mathbb{Q}$. This notation also applies to other spaces.
\begin{thm}[Existence (ii)]\label{existence2} If $(f, g, \xi)$ satisfies  {\rm\ref{as1}}, then there exists a solution of $(f, g, \xi)$ in $\mathscr{B}$.
\end{thm}
\begin{Proof}
\par
\noindent
(i).  We first show that it is equivalent to prove the existence result  given $|g_\cdot|\leq \frac{1}{8}$ $\mathbb{P}$-a.s. Suppose that $g$ is bounded by a positive constant $\tilde{g}$, that is, $|g_\cdot|\leq \tilde{g}$ $\mathbb{P}$-a.s.
Observe that,  for any $\theta >0$, $(Y, Z, N)$ is a solution  of $(f, g, \xi)$ if and only if $(\theta Y, \theta Z, \theta N)$  is a solution of $(f^\theta, g/\theta, \theta \xi)$, where
$f^\theta(t, y, z): = {\theta}f(t, \frac{y}{\theta},\frac{z}{\theta}).$
Obviously $f^\theta$ verifies \ref{as1} with the same Lipschitz coefficients as $f$. 
If we set $\theta  := 8  \tilde{g}$, then ${|g_\cdot/\theta|} \leq \frac{1}{8}$ $\mathbb{P}$-a.s. and hence satisfies the parametrization in Theorem \ref{existence1} (existence (i)). Therefore, we can and do assume  $|g_\cdot|\leq \frac{1}{8}$ $\mathbb{P}$-a.s. without loss of generality. 

(ii). Since $\norm{\xi}_\infty + \big\|{\big||f (\cdot, 0, 0)|\big|_T}\big\|_\infty < + \infty$, we can find  $n\in\mathbb{N}^+$ such that
\[
\xi=\sum_{i=1}^n \xi^i, \ 
f(t, 0, 0) = \sum_{i=1}^n f^i (t, 0, 0),
\]
where, for each $i \leq n$, $\xi^i$ is a  $\mathcal{F}_T$-measurable random variable, $f^i: \Omega \times [0, T]
\times \mathbb{R} \times \mathbb{R}^d \rightarrow \mathbb{R}$ is
$\Prog\otimes \mathcal{B}({\mathbb{R}})\otimes \mathcal{B}({\mathbb{R}^d})$-measurable and \[
\norm{\xi^i}^2_{\infty}+8 \big\|{\big||f^i(\cdot, 0, 0)|\big|_T}\big\|^2_{\infty} \leq \frac{1}{64}\exp \Big(       -\norm{A} 
\big(8\beta^2 \norm{{A}} + 8\gamma^2\big)
    \Big).
\]
Set $f^\prime  (t, y, z): = f(t, y, z)- f(t, 0, 0) $ and $(Y^0, Z^0\cdot M + N^0) \in \mathscr{B}$ such that $\norm{(Y^0, Z^0 \cdot M + N^0)} =0$. 
Now we use a recursion argument in the following way for $i=1, ..., n$.

By Theorem \ref{existence1}, there exists a solution $(Y^i, Z^i\cdot M + \widetilde{N}^{i}) \in \mathscr{B}(\mathbb{Q}^i)$ to the BSDE
\begin{align*}
Y_t^i &= \xi^i + \int_t^T \Big(f^i(s, 0, 0)+f^\prime(s, \sum_{j=0}^{i} Y_s^j, \sum_{j=0}^{i} Z_s^j) -  f^\prime(s, \sum_{j=0}^{i-1} Y_s^j, \sum_{j=0}^{i-1} Z_s^j) \Big)dA_s  \\ 
&+ \int_t^T g_s d \langle \widetilde{N}^{i}\rangle_s  -\int_t^T \big(Z_s^i dM_s +d\widetilde{N}_s^{i} \big),
\end{align*}
where 
\[
\frac{d\mathbb{Q}^i}{d\mathbb{P}}: = \mathcal{E}\Big(2g\cdot \sum_{j=0}^{i-1}N^j\Big)_T.
\]
Note that 
 the equivalent change of measure holds due to the fact that  $N^j \in \mathcal{M}^{BMO}$ for  $j \leq i-1$ and Theorem 2.3,  Kazamaki \cite{K1994}.
By Girsanov transformation and Theorem 3.6,  Kazamaki \cite{K1994},   $N^i := \widetilde{N}^{ i } + 2g \cdot \langle \widetilde{N}^{ i}, \sum_{j=0}^{i-1}N^j \rangle $ and $Z^i\cdot M$ belong to $\mathcal{M}^{BMO}$. This  further implies  $\langle N^i \rangle = \langle \widetilde{N}^i\rangle$ and 
$N^i = \widetilde{N}^{ i } + 2g \cdot \langle N^{ i}, \sum_{j=0}^{i-1}N^j \rangle .$ Hence $(Y^{i}, Z^{i}\cdot M + N^{i})\in \mathscr{B} $ solves
\begin{align*}
Y_t^i &= \xi^i + \int_t^T \Big(f^i (s, 0, 0 ) + f^\prime(s, \sum_{j=0}^{i} Y_s^j, \sum_{j=0}^{i} Z_s^j) -  f^\prime(s, \sum_{j=0}^{i-1} Y_s^j, \sum_{j=0}^{i-1} Z_s^j) \Big)dA_s  \\ 
&+ \int_t^T g_s d \Big(\langle N^i\rangle_s + 2\langle N^i, \sum_{j=0}^{i-1} N^j\rangle_s  \Big)  -\int_t^T \big(Z_s^i dM_s + dN_s^i\big ).
\end{align*}
Hence a recursion argument gives $(Y^i, Z^i, N^i)$ for $i =  1, ..., n$.

Define ${Y}: = \sum_{i=0}^n Y^i, Z:= \sum_{i=0}^n Z^i$ and $N:= \sum_{i=0}^n N^i$. 
 Clearly $(Y, Z\cdot M + N) \in \mathscr{B}$.
We show $(Y, Z, N)$ solves $(f, g, \xi)$.
In view of the definition of $f^\prime$, we 
sum up  the above BSDEs  to  obtain
\[
{Y}_t = \xi + \int_t^T \Big(\big(f(s, 0, 0) + f^\prime(s, {Y}_s, Z_s)\big)dA_s + g_s d\langle N \rangle_s\Big) -\int_t^T  \big( \delta Z_s dM_s + d\delta N_s\big).
\]
To conlcude the proof we use $f^\prime  (s, Y_s, Z_s): = f(s, Y_s, Z_s)- f(s, 0, 0) $.
\qed
\end{Proof}

We continue to show that  comparison theorem  and hence uniqueness also hold given Lipschitz-continuity. Similar results in different settings can be found, e.g., in 
\ \cite{MS2005}, \cite{HIM2005},  \cite{M2009}, \cite{T2008}.
\begin{thm}[Comparison]\label{compare1}
Let $(Y, Z\cdot M + N)$, $(Y^\prime, Z^\prime\cdot M+ N^\prime) \in \mathcal{S}^\infty \times \mathcal{M}^{BMO}$ be solutions of $(f, g, \xi)$, $(f^\prime, g^\prime, \xi^\prime)$, respectively. If 
$\mathbb{P}$-a.s. for any $(t, y, z) \in [0, T]\times\mathbb{R}\times\mathbb{R}^d$, 
$f(t, y ,z)\leq f^\prime(t, y, z)$, $g_t\leq g_t^\prime$, $\xi \leq \xi^\prime$ and $(f,g, \xi)$ verifies  {\rm\ref{as1}}, then $\mathbb{P}$-a.s. $Y_\cdot\leq Y^\prime_\cdot$.
\end{thm}
\begin{Proof}
 Set
$\delta Y := Y- Y^\prime$ and $\delta Z, \delta N, \delta \langle N\rangle, \delta \xi$, etc. analogously. For any $\tau\in\mathcal{T}$, $\mathbb{P}$-a.s. $f\leq f^\prime$ and $g_\cdot\leq g^\prime_\cdot$ imply by It\^{o}'s formula that 
\begin{align}
\delta Y_{t\wedge \tau} &= \delta Y_\tau + \int_{t\wedge \tau}^\tau \big(f(s, Y_s, Z_s)- f^\prime (s, Y_s^\prime, Z_s^\prime)\big)dA_s +\int_{t\wedge \tau}^\tau g_s d\langle N\rangle_s 
-\int_{t\wedge \tau}^\tau g^\prime_s d\langle N^\prime\rangle_s 
\nonumber\\
&-\int_{t\wedge \tau}^\tau \big(\delta Z_s dM_s + d\delta N_s\big ) \nonumber \\
& \leq  \delta Y_\tau + \int_{t\wedge \tau}^\tau\big( f(s, Y_s, Z_s) - f(s, Y^\prime_s, Z^\prime_s) \big)dA_s
  + \int_{t\wedge \tau}^\tau g^\prime_s d\delta \langle N\rangle_s  -\int_{t\wedge \tau}^\tau \big(\delta Z_s dM_s + d\delta N_s\big)\nonumber\\
&  =  \delta Y_\tau + \int_{t\wedge \tau}^\tau \big(\beta_s \delta Y_s + (\lambda_s \gamma_s)^\top (\lambda_s\delta Z_s) \big)dA_s 
  + \int_{t\wedge \tau}^\tau g^\prime_s d\delta \langle N\rangle_s  -\int_{t\wedge \tau}^\tau \big(\delta Z_s dM_s + d\delta N_s\big),\label{lqunique}
\end{align}
where  $\beta$ ($\mathbb{R}$-valued) and  $\gamma$ ($\mathbb{R}^d$-valued) are defined by
\begin{align*}
\beta_s &: =\mathbb{I}_{\{\delta Y_s \neq 0 \}} \frac{f(s, Y_s, Z_s) - f(s, Y_s^\prime, Z_s)}{\delta Y_s},\\
\gamma_s &: = \mathbb{I}_{\{\lambda_s\delta Z_s \neq \mathbf{0}\}}\frac{\big(f(s, Y_s^\prime, Z_s)-f(s, Y_s^\prime, Z_s^\prime)\big)\delta Z_s}{|\lambda_s \delta Z_s|^2},
\end{align*}
and  $\mathbf{0}:=(0, ..., 0)^\top.$ Note that $\gamma$ can be seen as defined in terms of discrete gradient. 
By \ref{as1}, $\beta_\cdot$ and $\int_0^\cdot \gamma^\top_s d\langle M \rangle_s \gamma_s$ are bounded processes,  hence $\gamma \cdot M \in \mathcal{M}^{BMO}.$ 
Given these facts we  use a change of measure to attain the comparison result. 
To this end, we define a BMO martingale
\[
\Lambda : = 
\gamma \cdot M  + g^\prime \cdot (N+N^\prime).
\]
In view of  Theorem 2.3 and Theorem 3.6, Kamazaki \cite{K1994},  we  define 
\[
\frac{d\mathbb{Q}}{d\mathbb{P}}: = \mathcal{E} ( \Lambda)_T.
\]
Hence
$ \delta N - g^\prime \cdot \delta \langle N \rangle$ and 
$\delta Z \cdot M- ( \gamma^\top\lambda^\top \lambda \delta Z) \cdot A $  belong to $\mathcal{M}^{BMO}(\mathbb{Q})$.
Therefore, (\ref{lqunique})  and $\mathbb{P}$-a.s. $\delta \xi \leq 0$ give
\begin{align*}
\delta Y_{t} &\leq \mathbb{E}^{\mathbb{Q}} \big[     \delta \xi\big|\mathcal{F}_{t} \big] + \mathbb{E}^{\mathbb{Q}}
\Big[ \int_{t}^{ T} \beta_s \delta Y_s dA_s \Big| \mathcal{F}_{t}\Big]  \\
&\leq     \mathbb{E}^{\mathbb{Q}}
\Big[ \int_{t}^{ T} \beta_s \delta Y_s dA_s \Big| \mathcal{F}_{t}\Big].
\end{align*}
Hence we obtain by Gronwall's lemma that  $\mathbb{P}$-a.s. $\delta Y_t \leq 0$.
Finally by the continuity of $Y$ and $Y^\prime$, we conclude that $\mathbb{P}$-a.s. $Y_\cdot \leq Y_\cdot^\prime.$
\qed
\end{Proof}

As a byproduct, we obtain the following existence and uniqueness result.
\begin{Corollary} [Uniqueness] \label{unique1}If $(f, g, \xi)$  satisfies  {\rm{\ref{as1}}}, then there exists a unique solution
in $\mathscr{B}$.
\end{Corollary}
\begin{Proof}
This is immediate from 
Theorem \ref{existence2} (existence (ii)) and 
 Theorem \ref{compare1} (comparison theorem). 
\qed
\end{Proof}
\subsection{Monotone Stability and Bounded Solutions of Quadratic BSDEs}
\label{section22}
In this section, we   prove a general  monotone stability result for quadratic BSDEs. 
Let us recall that Morlais \cite{M2009} uses  a stability-type argument for the existence result
   after performing an exponential transform which eliminates $g\cdot \langle N\rangle$.
But a  direct general stability result is not studied. Our work fills this gap. 

Secondly, as a byproduct of the stability property, we construct a bounded solution via 
  regularization through Lipschitz-quadratic BSDEs studied in Section \ref{section22}.  This  procedure is also called \emph{ Lipschitz-quadratic regularization} in the following context. Note that our definition of ``Lipschitz-quadratic'' is different from those in \cite{T2008}, \cite{BEK2013}.
To begin our proof, we give the assumptions for the whole section. 
\begin{as}\label{as2}
There exist $\beta \geq 0$, $\gamma >0$, an $\mathbb{R^+}$-valued $\Prog$-measurable process $\alpha$  and a continuous nondecreasing function $\varphi : \mathbb{R}^+ \rightarrow \mathbb{R}^+$ with $\varphi(0) = 0$ such that
${\norm{\xi}}_\infty  + \norm{|\alpha|_T}_\infty < +\infty$
and $\mathbb{P}$-a.s.
\begin{enumerate}
\item [(i)] for any $t\in[0, T]$, $(y, z)\longmapsto f(t, y, z)$ is continuous;
\item [(ii)]$f$ is monotonic at $y=0$, i.e., for any $(t, y, z)\in[0, T]\times\mathbb{R}\times\mathbb{R}^d$,
\[ \sgn (y) f(t, y, z) \leq \alpha_t + \alpha_t\beta |y|+ \frac{\gamma}{2}|\lambda_t z|^2;\]
\item [(iii)]  for any $(t, y, z)\in[0, T]\times\mathbb{R}\times\mathbb{R}^d$,
 \[ |f(t, y, z)| \leq \alpha_t + \alpha_t\varphi(|y|) + \frac{\gamma}{2}|\lambda_t z|^2. \]
\end{enumerate}
\end{as}

 We continue as before to call  $(\xi,|\alpha|_T)$ the \emph{data}.
 \ref{as2}(ii) 
  allows one to get rid of the linear growth  in $y$ which is required by  Kobylanski \cite{K2000} and Morlais \cite{M2009}. Assumption of this type for quadratic framework is motivated by Briand and Hu \cite{BH2008}.
Secondly, our results don't rely on the specific choice of $\varphi$.  Hence the growth condition in $y$  can be arbitrary as long as \ref{as2}(i)(ii)   hold. 

Given \ref{as2}, we first prove an a priori estimate. In order to treat $\langle Z\cdot M\rangle$ and $g \cdot \langle N \rangle$ more easily, we assume  $\mathbb{P}$-a.s. $|g_\cdot|\leq \frac{\gamma}{2}$ for the rest of this paper. 
\begin{Lemma}[A Priori Estimate]\label{aprioriestimate1} If $(f, g, \xi)$  satisfies  {\rm\ref{as2}} and $(Y, Z\cdot M + N) \in \mathcal{S}^\infty \times \mathcal{M}$  is a solution of $(f, g, \xi)$, then
\begin{align*}
\norm{Y} \leq \big\| e^{\beta |\alpha|_T}\big (|\xi| +|\alpha|_T\big)\big\|_\infty \end{align*}
and
\begin{align*}
&\norm {Z\cdot M + N }_{BMO} \leq c_{b},
\end{align*}
where $c_b$ is a constant only depending on $\beta, \gamma, \norm{\xi}_\infty, \norm{|\alpha|_T}_\infty$.
\end{Lemma}
\begin{Proof}  Set
 $u(x): = \frac{\exp({\gamma x})-1 - \gamma x}{\gamma^2}.$ The following auxiliary results will be useful:  
 $u(x) \geq 0, u^\prime(x) \geq 0$ and  $u^{\prime\prime}(x) \geq 1$ for $x\geq 0$;  $u(|\cdot|)\in \mathcal{C}^2(\mathbb{R})$ and $u^{\prime\prime}(x) = \gamma u^\prime(x) + 1$. 
For any  $\tau, \sigma \in \mathcal{T}$, It\^{o}'s formula yields
\begin{align*}
u(|Y_{\tau\wedge \sigma}|) =& u(|Y_{\sigma}|) + \int_{\tau \wedge \sigma}^\sigma u^\prime (|Y_s|)\sgn(Y_s)dY_s - \frac{1}{2}\int_{\tau\wedge \sigma}^\sigma u^{\prime\prime} (|Y_s|)\Big(Z_s^\top d\langle M\rangle_s Z_s + d\langle N\rangle_s\Big).
\end{align*}
By  \ref{as2}(ii),\[ 
\sgn(Y_s) f(s, Y_s, Z_s) \leq \alpha_s +\alpha_s \beta |Y_s|+ \frac{\gamma}{2}|\lambda_s Z_s|^2.\]
Note that $
\frac{\gamma}{2}u^\prime (|Y_s|)  -\frac{1}{2} u^{\prime\prime}(|Y_s|) = - \frac{1}{2}$,  $ 
{g_s}u^\prime (|Y_s|)  -\frac{1}{2} u^{\prime\prime}(|Y_s|) \leq - \frac{1}{2}.
$ and $u^\prime(|Y_s|) \leq \frac{e^{\gamma\norm{Y}}}{\gamma}$.
Hence, using these facts to the above equality yields
\begin{align*}
\frac{1}{2}\int_{\tau\wedge \sigma}^\sigma \Big( Z_s^\top d\langle M\rangle_s  Z_s + d\langle N \rangle_s \Big)&\leq 
\frac{e^{\gamma \norm{Y}}}{\gamma^2} +  \int_{\tau\wedge \sigma}^\sigma u^\prime (|Y_s|)\big(\alpha_s + \alpha_s \beta |Y_s|\big)dA_s \nonumber\\
& - \int_{\tau\wedge \sigma}^{ \sigma} u^{\prime}(|Y_s|)\sgn (Y_s) \big(Z_s dM_s + dN_s\big).  
\end{align*} 
To eliminate the local martingale,
we replace $\sigma$ by its localizing sequence and use Fatou's lemma to the left-hand side. Since  $Y^*$ and $|\alpha|_T$ are bounded random variables,  the right-hand side has a uniform constant  upper bound. Hence, we have 
\begin{align}
\frac{1}{2} \mathbb{E}\big[
\langle Z\cdot M  + N \rangle_{\tau, T} 
 \big|\mathcal{F}_\tau \big] \leq 
 \frac{e^{\gamma \norm{Y}}}{\gamma^2} + 
  \frac{e^{\gamma \norm{Y}}}{\gamma}(1+\beta \norm{Y}) \norm{|\alpha|_T}_\infty.
  \label{bmo}
\end{align}
 Now we turn to the estimate for $Y$. We fix $s\in [0, T]$ and for $t\in [s, T]$, set
\[
H_t : = \exp \Big(  \gamma e^{\beta |\alpha|_{s, t}}|Y_t| + \gamma \int_s^t e^{\beta |\alpha|_{s, u}}\alpha_u dA_u   \Big).
\]
We claim that $H$ is a submartingale. 
By Tanaka's formula,
\[
d|Y_t| = \sgn (Y_t) \big(Z_t dM_t + dN_t\big) -\sgn(Y_t)\big(f(t, Y_t, Z_t)dA_t + g_t d\langle N\rangle_t\big) + dL_t^0(Y),
\]
where $L^0(Y)$ is the local time of $Y$ at $0$. Hence, It\^{o}'s formula yields
\begin{align*}
dH_t &= \gamma H_t e^{\beta |\alpha|_{s, t}}\Big[  \sgn (Y_t)\big (Z_t dM_t + dN_t\big) \\
& + \Big(- \sgn(Y_t)f(t, Y_t, Z_t) +\alpha_t + \alpha_t \beta |Y_t| + \frac{\gamma}{2}e^{\beta |\alpha|_{s, t}}|\lambda_t Z_t|^2\Big)dA_t  \\
&+ \Big(- \sgn(Y_t)g_t +\frac{\gamma}{2}e^{\beta |\alpha|_{s, t}}\Big)d\langle N\rangle_t + dL_t^0(Y) 
\Big].
\end{align*}
By  \ref{as2}(ii) and  $ |g_\cdot|\leq \frac{\gamma}{2}$ again,  $(H_t)_{t\in[s, T]}$ is a bounded submartingale.  Hence, 
\begin{align*}
|Y_s| \leq \frac{1}{\gamma}\ln \mathbb{E} \big[H_T \big| \mathcal{F}_s\big].
\end{align*}
Thanks to the boundedness, we have 
\begin{align*}
\norm{Y} \leq \big\| e^{\beta |\alpha|_T} \big(|\xi| +|\alpha|_T\big)\big\|_\infty.
\end{align*}
Finally we  come back to (\ref{bmo}) 
and obtain the estimate for $Z\cdot M + N$.
\qed
\end{Proof}

Given the norm bound in Lemma \ref{aprioriestimate1}, 
we  turn to the main result of this section: the monotone stability result.  
 Later, as an immediate application,  we prove an existence result for quadratic BSDEs by Lipschitz-quadratic regularization. To start, we recall that $\mathcal{M}^2$ equipped with the norm $\norm{\widetilde{M}}_{\mathcal{M}^2}:=\mathbb{E}\big[\langle\widetilde{M}\rangle_T\big]^\frac{1}{2}$ for $\widetilde{M}\in\mathcal{M}^2$ is a Hilbert space. 
\begin{thm}[Monotone Stability] \label{monotonestability}
Let $(f^n, g^n, \xi^n)_{n\in\mathbb{N}^+}$  satisfy  {\rm\ref{as2}} associated with $(\alpha, \beta, \gamma, \varphi)$, and  $(Y^n, Z^n \cdot M + N^n)$ be their solutions in $\mathscr{B}$, respectively.  Assume
\begin{enumerate}
\item [{\rm(i)}]
 $Y^n$  is monotonic in $n$ and  $
\xi^n - \xi \longrightarrow 0           \ \mathbb{P}$-a.s. with $\sup_n \norm{\xi^n}_\infty< +\infty${\rm;} 
\item [{\rm(ii)}]$\mathbb{P}$-a.s. for any $t\in [0, T]$, $g^n_t - g_t \longrightarrow 0${\rm;}  
\item [{\rm(iii)}] $\mathbb{P}$-a.s. for any $t\in[0, T]$ and $y^n \longrightarrow y, z^n \longrightarrow z$,
$f^n(t, y^n, z^n)\longrightarrow f(t, y, z)$.
\end{enumerate}

Then there exists a process $(Y, Z\cdot M + N) \in \mathscr{B}$ such that $Y^n$ converges to $Y$ $\mathbb{P}$-a.s. uniformly on $[0, T]$ and $(Z^n \cdot M + N^n)$ converges to $(Z\cdot M + N)$ in  $\mathcal{M}^2$ as $n$ goes to $+\infty$. Moreover, $(Y, Z, N)$ solves $(f, g, \xi)$.

\end{thm}
\begin{Proof}
Without loss of generality we only consider  $Y^n$ to be increasing in $n$.
By Lemma \ref{aprioriestimate1} (a priori estimate),
\begin{align}
\sup_n  \norm{Y^n}+\sup_n \norm{Z^n \cdot M + N^n}_{BMO}  \leq c_b, \label{e31}
\end{align}
where $c_b$ is a constant  only depending on $\beta, \gamma, \sup_n{\norm{\xi^n}}_\infty, \norm{|\alpha|_T}_\infty$. We  rely intensively on the boundedness  result in (\ref{e31})  to derive the limit.

(i). We prove  the  convergence of the solution sequences. 
 Due to (\ref{e31}),
there exists a bounded monotone limit 
$
Y_t:= \lim_{n} Y_t^n,
$   a subsequence 
 indexed by $\{n_k\}_{k\in\mathbb{N}^+} \subseteq \mathbb{N}^+$ 
 and $Z\cdot  M +N \in \mathcal{M}^2$
 such that 
 $Z^{n_k}\cdot M + N^{n_k}$ converges weakly 
in $\mathcal{M}^2$
 to $Z\cdot M + N$
as $k$ goes to $+\infty.$ 
 The remaining task is to show $Z\cdot M + N$  is the $\mathcal{M}^2$-limit of the whole sequence.  To this end, we  define
$
u(x) := \frac{\exp({8\gamma x})-8\gamma x -1}{64\gamma^2}.
$
Recall that $u(x) \geq 0, u^\prime (x)\geq 0$ and $u^{\prime\prime}(x) \geq 0$ for $x\geq 0$; $u \in \mathcal{C}^2(\mathbb{R})$ and $u^{\prime\prime}(x)= 8\gamma u^\prime(x)+1$.  For any $  m\in \{n_k\}_{k\in\mathbb{N}^+}$ ,$n \in \mathbb{N}^+$ with $m \geq n$, define $\delta Y^{m, n}: = Y^m - Y^n, \delta Y^{n}: = Y - Y^n$ and $\delta Z^{m, n}, \delta Z^n, \delta N^{m, n}, \delta N^n$, etc. analogously.
By It\^{o}'s formula,  we have
\begin{align}
\mathbb{E}\big[u(\delta Y_0^{m, n})\big]- \mathbb{E}\big[u(\delta \xi^{m, n})\big]
&=\mathbb{E}\Big[\int_0^T  u^\prime(\delta Y_s^{m, n})\big(f^m (s, Y^m_s, Z^m_s)- f^n (s, Y^n_s, Z^n_s)\big)dA_s\Big]
\nonumber\\
&+ \mathbb{E}\Big[\int_0^T u^\prime (\delta Y_s^{m, n}) \big(g_s^m d\langle N^m \rangle_s  
-  g_s^{n} d \langle N^n \rangle_s \big)\Big]\nonumber \\
&-\frac{1}{2}\mathbb{E}\Big[  \int_0^T u^{\prime\prime}(\delta Y_s^{m, n} )  \Big( (\delta Z_s^{m, n})^\top d\langle M\rangle_s (\delta Z_s^{m, n})+ d\langle \delta N^{m, n}\rangle_s \Big) \Big]. \label{e32}
\end{align}
Since $f^m$ and $f^{n}$ verify \ref{as2} associated with $(\alpha, \beta, \gamma, \varphi)$, we have
\begin{align*}
|f^m (s, Y^m_s, Z^m_s) &- f^n (s, Y^n_s, Z^n_s)|  \nonumber\\
&\leq \alpha_s^\prime +\frac{\gamma}{2}|\lambda_s Z^m_s|^2 + \frac{\gamma}{2}|\lambda_s Z^n_s|^2 \nonumber\\
& \leq \alpha_s^\prime+ \frac{3\gamma}{2}\big(|\lambda_s \delta Z_s^{m, n}|^2 + |\lambda_s \delta Z^n_s|^2 + |\lambda_s Z_s|^2\big) + \gamma \big(|\lambda_s \delta Z_s^n|^2 + |\lambda_s Z_s|^2\big) \nonumber \\
& \leq \alpha_s^\prime +\frac{3\gamma}{2}|\lambda_s \delta Z^{m, n}_s|^2 +  \frac{5\gamma}{2} \big( |\lambda_s \delta Z^n_s |^2 + |\lambda_s Z_s|^2 \big),
\end{align*}
where
\[
\alpha_s^\prime :  = 2\alpha_s\big(1+\varphi(c_b) \big) \geq   2\alpha_s + \alpha_s\varphi(|Y^n_s|)+ \alpha_s\varphi(|Y^m_s|).
\]
Moreover,
\begin{align*}
g^m d\langle N^m \rangle - g^n d\langle N^n \rangle  &\ll
\frac{\gamma}{2}d\langle N^m \rangle + \frac{\gamma}{2} d\langle N^n \rangle \\
&\ll
\frac{3\gamma}{2}d\langle \delta N^{m , n}\rangle + \frac{5\gamma}{2}\big(d\langle \delta N^n \rangle
+ d\langle N \rangle  \big).
\end{align*}
Plugging the above inequalities into (\ref{e32}), we deduce that 
\begin{align}
&\mathbb{E}\Big[ \int_0^T  \Big(\frac{1}{2}u^{\prime\prime} -\frac{3\gamma}{2}u^\prime\Big) (\delta Y_s^{m, n})  |\lambda_s \delta Z_s^{m, n}|^2 dA_s	\Big]
+ \mathbb{E}\Big[  \int_0^T \Big(\frac{1}{2} u^{\prime\prime}-\frac{3\gamma}{2}u^\prime\Big)(\delta Y_s^{m, n})   d\langle\delta N^{m, n}\rangle_s     \Big]\nonumber \\
&\leq
\mathbb{E}\big[u(\delta \xi^{m, n})\big] + \mathbb{E}\Big[\int_0^T u^\prime (\delta Y^{m, n}_s) \Big( \alpha_s^\prime +\frac{5\gamma}{2}\big(|\lambda_s \delta Z_s^{n}|^2  + |\lambda_s Z_s|^2 \big)  \Big)  dA_s \Big]\nonumber\\
&+  \mathbb{E}\Big[\int_0^T u^\prime (\delta Y_s^{m, n}) \frac{5\gamma}{2} \big(d\langle \delta N^n \rangle_s + d\langle N \rangle_s    \big) \Big]
 \label{e34} 
\end{align}
Due to the weak convergence result and  convexity of $z\longmapsto |z|^2$, $N\longmapsto \langle N\rangle$, we obtain
\begin{align*}
\mathbb{E}\Big[ \int_0^T  \Big(\frac{1}{2}u^{\prime\prime} -\frac{3\gamma}{2}u^\prime\Big) (\delta Y_s^{ n})  |\lambda_t Z_s^{ n}|^2 dA_s	\Big] &\leq \liminf_{m} 
\mathbb{E}\Big[ \int_0^T  \Big(\frac{1}{2}u^{\prime\prime} -\frac{3\gamma}{2}u^\prime\Big) (\delta Y_s^{m, n})  |\lambda_t Z_s^{ m, n}|^2 dA_s	\Big],\\
 \mathbb{E}\Big[  \int_0^T \Big(\frac{1}{2}u^{\prime\prime}-\frac{3\gamma}{2}u^\prime \Big)(\delta Y_s^{ n})   d\langle\delta N^{ n}\rangle_s     \Big]&\leq
 \liminf_m\mathbb{E}\Big[  \int_0^T \Big(\frac{1}{2}u^{\prime\prime}-\frac{3\gamma}{2}u^\prime\Big)(\delta Y_s^{m, n})   d\langle\delta N^{m, n}\rangle_s     \Big].
\end{align*}
We then come back to 
 (\ref{e34}) and send $m$ to  $+\infty$ along $\{n_k \}_{k\in \mathbb{N}^+}$. Taking the  above inequalities into account and using $u^{\prime}(\delta Y^{m, n}_s)\leq u^{\prime}(\delta Y^{n}_s)$
 to the right-hand side, 
 (\ref{e34}) becomes
\begin{align}
&\mathbb{E}\Big[ \int_0^T  \Big(\frac{1}{2}u^{\prime\prime} -\frac{3\gamma}{2}u^\prime\Big) (\delta Y_s^{ n})  |\lambda_s Z_s^{ n}|^2 dA_s	\Big]\nonumber\\
&+ \mathbb{E}\Big[\int_0^T \Big(\frac{1}{2}u^{\prime\prime}-\frac{3\gamma}{2}u^{\prime}\Big)(\delta Y_s^n) d\langle \delta N^n \rangle_s  \Big]
\nonumber\\&\leq \mathbb{E}\big[u(\delta \xi^{ n})\big]  + \mathbb{E}\Big[\int_0^T u^\prime (\delta Y^{ n}_s) \Big( \alpha^\prime_s +\frac{5\gamma}{2}\big(|\lambda_s \delta Z_s^{n}|^2  + |\lambda_s Z_s|^2 \big)  \Big)  dA_s \Big]\nonumber\\
&+ \frac{5\gamma}{2}\mathbb{E}\Big[\int_0^T u^\prime (\delta Y_s^{n})  \big(d\langle \delta N^n \rangle_s + d\langle N \rangle_s    \big) \Big].
\label{e35}
\end{align}
Since $u^{\prime\prime} (x) - 8\gamma u^\prime (x) =1$,   rearranging terms give
\begin{align}
&\frac{1}{2}E\big[ \big(\delta N_T^n\big)^2 \big]  + \frac{1}{2}\mathbb{E}\Big[ \int_0^T |\lambda_s \delta Z_s^n|^2 dA_s\Big]
 \nonumber\\
&\leq \mathbb{E}\big[u(\delta \xi^{ n})\big]+ \mathbb{E}\Big[\int_0^T u^\prime (\delta Y^{ n}_s) \Big( \alpha^\prime_s +\frac{5\gamma}{2} |\lambda_s Z_s|^2   \Big)  dA_s \Big]+ \frac{5\gamma}{2}\mathbb{E}\Big[ \int_0^T u^{\prime}(\delta Y_s^n)d\langle N \rangle_s\Big].
 \label{e36}
\end{align}
Finally, by sending $n$ to $+\infty$ and dominated convergence we deduce the convergence.

(ii).
We  prove $(Y, Z\cdot M + N)\in \mathscr{B}$ and solves $(f, g, \xi)$.
Here we rely on the same arguments as in  Kobylanski \cite{K2000} or  Morlais \cite{M2009} and omit the details here.
In addition to their deductions, we need to prove the $u.c.p$ convergence of  $g^n\cdot \langle N^n \rangle$, which holds if 
\begin{align*}
\lim_{n\rightarrow \infty} \mathbb{E}\Big[ \Big| \int_0^\cdot \big(g^n_s d\langle N^n \rangle_s - g_s d\langle N\rangle_s \big)\Big|^* \Big] =0.
\end{align*}
Indeed, by Kunita-Watanabe inequality and Cauchy-Schwartz inequality, 
\begin{align*}
 \mathbb{E}\Big[ \Big| \int_0^\cdot \big(g^n_s d\langle N^n \rangle_s &- g_s d\langle N \rangle_s \big) \Big|^* \Big]
= \mathbb{E}\Big[ \Big| \int_0^\cdot \Big(g^n_s d\big(\langle N^n \rangle_s - \langle N \rangle_s \big)+ ( g^n_s-g_s )d\langle N \rangle_s \Big) \Big|^* \Big]\\
&\leq \frac{\gamma}{2} \mathbb{E}\big[     \langle N^n - N\rangle_T  \big]^{\frac{1}{2}} \mathbb{E}\big[     \langle N^n + N\rangle_T  \big]^{\frac{1}{2}} + \mathbb{E}\Big[ \Big|  \int_0^\cdot (g^n_s -g_s)d\langle N\rangle_s \Big|^*   \Big]\\
&\leq  \gamma c_b \mathbb{E}\big[     \langle N^n - N\rangle_T  \big]^{\frac{1}{2}} + \mathbb{E}\Big[  \int_0^T |g^n_s -g_s|d\langle N\rangle_s   \Big].
\end{align*}
We then conclude by  $\mathcal{M}^2$-convergence of $N^n$ and dominated convergence used to the second term. Finally 
$Z\cdot M + N \in\mathcal{M}^{BMO}$ by Lemma \ref{aprioriestimate1} (a priori estimate).

For decreasing $Y^n$, we take $m \in \mathbb{N}^+, n \in \{n_k\}_{k\in \mathbb{N}^+}$ with $n \geq m$ and conclude with exactly the same arguments.
\qed
\end{Proof}

There are several major improvements compared to existing monotone stability results. First of all, in contrast to  Kobylanski \cite{K2000} and  Morlais \cite{M2009}, we get rid of linear growth in $y$ by merely assuming \ref{as2}, and allow $g$ to  be any  bounded process.
Secondly, we treat the convergence in a more direct and general way than Morlais \cite{M2009}.

Another advantage concerns the existence result. 
 Thanks to Section \ref{section21} and Theorem \ref{monotonestability},
 we are able to perform directly a Lipschitz-quadratic regularization without exponential transforms; this is in contrast to Morlais \cite{M2009}.
One can also benefit from our stability result in obtaining the existence results for unbounded solutions with more flexible assumptions;  see Section \ref{section23}.

\begin{prop}[Existence]\label{existence3}
If $(f, g, \xi)$ satisfy  {\rm\ref{as2}}, then there exists a solution in $\mathscr{B}$.
\end{prop}
\begin{Proof}
We use a double approximation procedure and use Theorem \ref{monotonestability}  (monotone stability) to take the limit. 
 Define
\begin{align*}
f^{n, k}(t, y, z): &= \inf_{y^\prime, z^\prime} \big\{f^+ (t, y^\prime, z^\prime) + n|y-y^\prime| + n|\lambda_t(z-z^\prime)|\big \}\\
&- \inf_{y^\prime, z^\prime}\big\{ f^- (t, y^\prime, z^\prime) + k|y-y^\prime| + k|\lambda_t(z-z^\prime)|\big\}.
\end{align*}
By
Lepeltier and  San Martin \cite{LS1997}, $f^{n, k}$ is Lipschitz-continuous in $(y, z)$; as $k$  goes to $+\infty$, $f^{n, k}$ converges increasingly   uniformly on compact sets to a limit denoted by $f^{n, \infty}$;  as $n$  goes to $+\infty$, $f^{n, \infty}$ converges increasingly uniformly on compact sets to $f$.

By Corollary \ref{unique1}, there exists a unique solution $(Y^{n ,k}, Z^{n, k}\cdot M + N^{n, k})\in\mathscr{B}$ to $(f^{n, k}, g, \xi)$; by  Theorem \ref{compare1} (comparison theorem), $Y^{n, k}$ is increasing in $n$ and decreasing in $k$, and is uniformly bounded due to
Lemma \ref{aprioriestimate1} (a priori estimate).  We then fix $n$ and use  Theorem \ref{monotonestability} to the sequence indexed by $k$ to obtain a solution
 $(Y^{n},Z^{n}\cdot M + N^{n} ) \in$ $\mathscr{B}$ to $(f^{n, \infty}, g, \xi)$. Due to the $\mathbb{P}$-a.s. uniform convergence of $Y^{n, k}$ we  can pass the comparison property to $Y^{n}$. We use Theorem \ref{monotonestability} again to conclude. 
\qed
\end{Proof}
\begin{remm}
In contrast to  Kobylanski \cite{K2000}, the existence of a maximal or minimal solution is not available (yet) given \ref{as1} as the double approximation procedure makes the comparison between solutions impossible.
\end{remm}

There is also a rich literature on the uniqueness of a bounded solution of quadratic BSDEs; see, e.g.,  \cite{K2000}, \cite{MS2005}, \cite{HIM2005},  \cite{M2009}.  Roughly speaking,
they essentially  rely
a type of locally Lipschitz-continuity and use a change of measure  analogously to Section \ref{section21}. 
 The proof in our setting is exactly the same and hence  omitted to save pages. 

To end this section, we briefly present various structure conditions used in different situations.

\begin{asp}{as2} 
 \label{as2'} 
 There exist $\beta \geq 0, \gamma >0$, an $\mathbb{R}^+$-valued $\Prog$-measurable process $\alpha$, 
 and a continuous nondecreasing function $\varphi : \mathbb{R}^+ \rightarrow \mathbb{R}^+$ with $\varphi(0) = 0$ such that
$\mathbb{P}$-a.s.  
\begin{enumerate}
\item [(i)] for any $t\in [0, T]$, $(y, z)\longmapsto f(t, y, z)$ is continuous;
\item [(ii)] $f$ is monotonic  at $y=0$, i.e., for any $(t, y, z)\in[0, T]\times\mathbb{R}\times\mathbb{R}^d$,
\begin{align*}\sgn (y) f(t, y, z) \leq \alpha_t + \beta |y|+ \frac{\gamma}{2}|\lambda_t z|^2;\end{align*}
\item [(iii)]  for any $(t, y, z)\in[0, T]\times\mathbb{R}\times\mathbb{R}^d$,
\[
|f(t, y, z)| \leq \alpha_t + \varphi(|y|)+ \frac{\gamma}{2}|\lambda_t z|^2. \] 
\end{enumerate}
\end{asp}

Given bounded data, \ref{as2'} implies \ref{as2}. Indeed, 
\begin{align*}
\sgn(y) f(t, y, z )&\leq \alpha_t \vee 1 +  (\alpha_t \vee 1)\beta |y| + \frac{\gamma}{2}|\lambda_t z|^2, \\
|f(t, y, z)|  &\leq \alpha_t \vee 1 + (\alpha_t\vee 1)\varphi(|y|)+ \frac{\gamma}{2}|\lambda_t z|^2.
\end{align*}
Hence \ref{as2'} verifies  \ref{as2} associated with $(\alpha \vee 1, \beta, \gamma, \varphi)$. 
However, given unbounded data,   \ref{as2'}
appears to be more natural and convenient.
 This will be discussed in detail in Section \ref{section23}.

In  particular situations where the estimate for $\int_0^T |f(s, Y_s, Z_s)|dA_s$ is needed,  e.g., in analysis of measure change (see Section \ref{section24}), there has to be a linear growth in $y$, which corresponds to the following assumption 
\begin{aspp}{as2}
\label{as2''}
There exist  $\beta \geq 0$, $\gamma >0$, an $\mathbb{R}^+$-valued $\Prog$-measurable process $\alpha$ such that $\mathbb{P}$-a.s.
\begin{enumerate}
\item [(i)] for any $t\in [0, T]$, $(y, z)\longmapsto f(t, y, z)$ is continuous;
\item [(ii)]
for any $(t, y, z)\in[0, T]\times\mathbb{R}\times\mathbb{R}^d$,
\[
|f(t, y, z)| \leq \alpha_t + \beta|y| + \frac{1}{2}|\lambda_t z|^2.
\]
\end{enumerate}
Indeed, \ref{as2''}
 enables one to obtain the estimate for
$
\int_0^T |f(s, Y_s, Z_s)|dA_s
$
via
\[
\int_0^T |f(s, Y_s, Z_s)|dA_s 
\leq |\alpha|_T + \beta \norm{A}Y^*
+ \frac{\gamma}{2}\langle Z\cdot M \rangle_T.
\]
\end{aspp}

\subsection{Unbounded Solutions of Quadratic BSDEs}
\label{section23}
This  section extends  Section \ref{section21}, \ref{section22} to unbounded solutions. We prove an existence result and later  show that the uniqueness holds given convexity assumption as an additional requirement.
We point out that similar results have been obtained by Mocha and Westray \cite{MW2012}, but our results rely on much fewer  assumptions and are more natural.
Analogously to section \ref{section22}, we give an a priori estimate in the first step.
We keep in mind that  
$\mathbb{P}$-a.s. $ |g_\cdot| \leq \frac{\gamma}{2}$ throughout our study.
\begin{Lemma}[A priori estimate]\label{aprioriestimate2} If $(f, g, \xi)$ satisfies  {\rm\ref{as2'}} and
$(Y, Z\cdot M + N)\in \mathcal{S}\times\mathcal{M}$ is a solution of $(f, g, \xi)$ such that the process 
\[
\exp \Big( \gamma e^{\beta A_T}{ |Y_\cdot|} + \gamma \int_0^T e^{\beta A_s}\alpha_s dA_s    \Big)
\]
is of class $\mathcal{D}$, then 
\begin{align}
|Y_s| \leq \frac{1}{\gamma} \ln \mathbb{E}\Big[ \exp\Big(\gamma e^{\beta A_{s, T}}|\xi| + 
\gamma \int_s^T e^{\beta A_{s, u}}{\alpha_u dA_u} \Big)\Big| \mathcal{F}_s \Big]. \label{ey}
\end{align}
\end{Lemma}
\begin{Proof}
We  fix $s\in [0, T]$, and for $t\in [s, T]$, set
\begin{align}
H_t : = \exp \Big(  \gamma e^{\beta A_{s, t}}|Y_t| + \gamma \int_s^t e^{\beta A_{s, u}}\alpha_u dA_u   \Big). \label{H}
\end{align}
We claim that $H$ is a local submartingale.  Indeed, 
by Tanaka's formula
\[
d|Y_t| = \sgn (Y_t) \big(Z_t dM_t + dN_t \big) -\sgn(Y_t)\big(f(t, Y_t, Z_t)dA_t + g_t d\langle N\rangle_t\big) + dL_t^0(Y),
\]
where $L^0(Y)$ is the local time of $Y$ at $0$. Hence, It\^{o}'s formula yields
\begin{align*}
dH_t &= \gamma H_t e^{\beta A_{s, t} }\Big[  \sgn (Y_t) \big(Z_t dM_t + dN_t \big)  \\
& + \Big(- \sgn(Y_t)f(t, Y_t, Z_t) +\alpha_t +  \beta |Y_t| + \frac{\gamma}{2}e^{\beta  A_{s, t} }|\lambda_t Z_t|^2\Big)dA_t  \\
&+ \Big(- \sgn(Y_t)g_t +\frac{\gamma}{2}e^{\beta A_{s, t} }\Big)d\langle N\rangle_t + dL_t^0(Y) 
\Big].
\end{align*}
By  \ref{as2'}(ii), $H$ is a local submartingale. To eliminate the local martingale part, we replace $\tau$ by its localizing sequence on $[s, T]$, denoted by
 $\{\tau_n\}_{n\in \mathbb{N}^+}$. Therefore, 
\begin{align*}
|Y_s| &\leq \frac{1}{\gamma}\ln \mathbb{E} \big[H_{T\wedge \tau_n} \big| \mathcal{F}_s\big]\\
&\leq \frac{1}{\gamma} \ln \mathbb{E}\Big[ \exp\Big(\gamma e^{\beta A_{s, T\wedge \tau_n}}|Y_{T\wedge \tau_n}| + 
\gamma \int_s^{T\wedge \tau_n} e^{\beta A_{s, u}}{\alpha_u dA_u} \Big)\Big| \mathcal{F}_s \Big]. 
\end{align*}
Finally by class $\mathcal{D}$ property
we conclude by sending $n$ to $+\infty$.
\qed
\end{Proof}

We then know from Lemma \ref{aprioriestimate2} that exponential moments integrability on $|\xi| + |\alpha|_T$ is a natural requirement for the existence result. 
\begin{remm}
 \ref{as2'} addresses the issue of integrability better than
  \ref{as2}.  
To show this, let us assume \ref{as2}. We then deduce from Lemma \ref{aprioriestimate1}  and corresponding class $\mathcal{D}$ property that 
\begin{align}
|Y_s| \leq \frac{1}{\gamma} \ln \mathbb{E}\Big[  \exp\Big(\gamma e^{\beta|\alpha|_{s, T}}|\xi| + \gamma \int_s^T e^{\beta |\alpha|_{s, u}} \alpha_u dA_u     \Big)         \Big| \mathcal{F}_s     \Big]. \label{bound2}
\end{align}
Obviously, in (\ref{bound2}), even exponential moments integrability  
is not sufficient to ensure the well-posedness of the a priori estimate.
For more dicusssions on the choice of structure conditions, the reader shall refer to  Mocha and Westray \cite{MW2012}.
\end{remm}

Motivated by the above discussions,  we prove  an existence  result given \ref{as2'} and exponential moments integrability. Analogously to Theorem \ref{existence3}, 
we  use a Lipschitz-quadratic regularization 
and take the limit by the monotone stability result in Section \ref{section22}.
The a priori bound for $Y$ obtained in  Lemma \ref{aprioriestimate2}  is also  crucial to the construction of an unbounded solution. 
 \begin{thm}[Existence]
\label{existence4}
 If $(f, g, \xi)$ satisfies {\rm\ref{as2'}} and $e^{\beta A_T} \big(|\xi|+ |\alpha|_T\big)$ has exponential moment of order $\gamma$, i.e.,  
\[
\mathbb{E}\Big[ \exp\Big( \gamma e^{\beta A_T}\big (|\xi|+ |\alpha|_T \big)      \Big)\Big] < + \infty,
\]
then there exists a solution  verifying {\rm(\ref{ey})}.
\end{thm}
\begin{Proof} We introduce the notations used throughout the proof.
Define the process \[
X_t: = \frac{1}{\gamma}\ln \mathbb{E}\Big[ \exp \Big(\gamma e^{\beta A_T}\big( |\xi| + |\alpha|_T\big ) \Big)\Big| \mathcal{F}_t  \Big].
\]
Obviously $X$  is continuous by the  continuity of the filtration. 
For $m, n \in \mathbb{N}^+$, set 
\begin{align*}
\tau_m &:= \inf \big\{ t\geq 0: |\alpha|_t + X_t \geq m \big\} \wedge T,\\
\sigma_n &: = \inf \big\{ t\geq 0: |\alpha|_t \geq n\big\} \wedge T.
\end{align*}
It then follows from the  continuity of $X$ and  $|\alpha|_\cdot$ that $\tau_m$ and $\sigma_n$ increase  stationarily to $T$ as $m, n$ goes to $+\infty$, respectively.
To  apply a double approximation procedure, we  define 
 \begin{align*}
f^{n, k}(t, y, z)&: =\mathbb{I}_{\{t\leq \sigma_n \}} \inf_{y^\prime, z^\prime}\big \{ f^+(t, y^\prime, z^\prime) + n|y-y^\prime|+n|\lambda_t(z-z^\prime)| \big\} \\
&-\mathbb{I}_{\{t\leq \sigma_k \}} \inf_{y^\prime, z^\prime} \big\{ f^-(t, y^\prime, z^\prime) + k|y-y^\prime|+k|\lambda_t(z-z^\prime)| \big\},
\end{align*}
and $
\xi^{n , k}: = \xi^+ \wedge n - \xi^- \wedge k.$   

Before proceeding to the proof we give some useful facts. 
By  Lepeltier and San Martin \cite{LS1997}, 
$f^{n, k}$ is  Lipschitz-continuous in $(y, z)$; as  $k$ goes to $+\infty$,  $f^{n, k}$ converges decreasingly 
 uniformly on compact sets to a limit denoted by $f^{n,\infty}$; as $n$ goes to $+\infty$, $f^{n, \infty}$ converges increasingly uniformly on compact sets to $F$. Moreover, $\big||f^{n, k}(\cdot, 0, 0) |\big|_T$ and $\xi^{n, k}$ are bounded. 
 
Hence, by Corollary \ref{unique1},  there exists a unique solution $(Y^{n,k}, Z^{n, k}\cdot M + N^{n, k})\in\mathscr{B}$ to $(f^{n, k}, g, \xi^{n, k})$; by Theorem \ref{compare1} (comparison theorem), $Y^{n, k}$ is increasing in $n$ and decreasing in $k$.
Analogously to Proposition  \ref{existence3}, we wish to take the limit by  
 Theorem \ref{monotonestability} (monotone stability). 
 
 However, 
 $|f^{n, k}(\cdot, 0, 0)|_T$ and $\xi^{n, k}$ 
are not uniformly bounded in general.
To overcome this difficulty, we use Lemma \ref{aprioriestimate2} (a priori estimate) and work on random interval where	
$Y^{n,k}$ and $|f^{n, k}(\cdot, 0, 0)|_\cdot$ are uniformly bounded. This is the motivation to introduce $X$ and $\tau_m$. To be more precise, the localization procedure is 
as follows.

Note that   $(f^{n, k}, g, \xi^{n, k})$ verifies  \ref{as2'} associated with $(\alpha, \beta, \gamma, \varphi)$.  $Y^{n, k}$ being bounded implies that it is of class $\mathcal{D}$. 
Hence  from Lemma \ref{aprioriestimate2}  we have 
\begin{align}
|Y^{n, k}_t|&  \leq \frac{1}{\gamma}\ln \mathbb{E}\Big[ \exp \Big(\gamma e^{\beta A_{t, T}} |\xi^{n ,k}| 
+\gamma \int_{t}^{T} e^{\beta A_{t, s} } \alpha_{s} \mathbb{I}_{\{ s\leq \sigma_n \wedge \sigma_k\}}dA_s
 \Big)\Big| \mathcal{F}_t  \Big]
\nonumber\\
&
\leq \frac{1}{\gamma} \ln \mathbb{E}\Big[ \exp\Big(\gamma e^{\beta A_{t, T}}|\xi| + 
\gamma \int_t^T e^{\beta A_{t, T}}{\alpha_s dA_s} \Big)\Big| \mathcal{F}_t \Big] \label{x1} \\
&\leq X_t.  \nonumber
\end{align}
In view of the definition of $\tau_m$, we have 
\begin{align}
|Y^{n, k}_{t\wedge \tau_m} |&\leq X_{t\wedge \tau_m} \leq m, \nonumber \\
\big||f^{n, k}(\cdot, 0, 0)|\big|_{\tau_m} &\leq |\mathbb{I}_{[0, \tau_m]}\alpha|_{\tau_m} \leq m.
 \label{bound}
\end{align}
Hence $\big||f^{n, k}(\cdot, 0, 0) |\big|_\cdot$ and $Y^{n, k}$ are uniformly bounded on $[0, \tau_m]$. Secondly, given $(Y^{n, k}, Z^{n, k}\cdot M + N^{n, k})$ which solves $(f^{n, k}, g, \xi^{n, k})$, it is immediate that 
$(Y^{n, k}_{\cdot \wedge \tau_m}, (Z^{n, k}\cdot M + N^{n, k})_{\cdot \wedge \tau_m})$  solves $(\mathbb{I}_{[0, \tau_m ]}(t)f^{n, k}(t, y, z), g, Y_{\tau_m}^{n, k})$. We then use Theorem \ref{monotonestability}  as in Proposition \ref{existence3} to construct a pair $(\widetilde{Y}^{m}, (\widetilde{Z}^m \cdot M + \widetilde{N}^m))$ which solves 
$(f, g, \sup_n\inf_k Y_{\tau_m}^{n, k})$, i.e., 
\begin{align}
\widetilde{Y}^m_t =\sup_n\inf_k Y_{\tau_m}^{n, k} + \int_{t\wedge \tau_m}^{\tau_m} \big(
F(s,\widetilde{Y}^m_s, \widetilde{Z}^m_s ) dA_s + g_s \langle \widetilde{N}^m \rangle_s \big)
 -\int_{t\wedge \tau_m}^{\tau_m}   \big(\widetilde{Z}_s^m dM_s + d\widetilde{N}_s \big). \label{bsdem}
\end{align}
Moreover, 
  $\widetilde{Y}^m$ is the $\mathbb{P}$-a.s. uniform limit of $Y^{n, k}_{\cdot\wedge \tau_m}$ and 
$\widetilde{Z}^m\cdot M + \widetilde{N}^m$ is the $\mathcal{M}^2$-limit of $(Z^{n, k}\cdot M + N^{n, k})_{\cdot \wedge \tau_m}$ as $k, n$ go to $+\infty$. Hence
\begin{align}
\widetilde{Y}^{m+1}_{\cdot\wedge \tau_m}&= \widetilde{Y}^{m}_{\cdot\wedge \tau_m}\  \mathbb{P}\text{-a.s.}, \nonumber\\
  \mathbb{I}_{\{t\leq \tau_m\}}\lambda_t\widetilde{Z}^{m+1}_t &= \lambda_t\widetilde{Z}^m_t\ dA \otimes d\mathbb{P}\text{-a.e},\nonumber
  \\
  \widetilde{N}^{m+1}_{\cdot\wedge \tau_m} &= \widetilde{N}^m_{\cdot\wedge \tau_m}\  \mathbb{P}\text{-a.s.} \label{ae}
\end{align}
Define $(Y, Z, N)$ on $[0, T]$ by
\begin{align*}
Y_t &:= \mathbb{I}_{\{t\leq \tau_m\}} \widetilde{Y}_t^1  + \sum_{m\geq 2} \mathbb{I}_{]\tau_{m-1}, \tau_m]}\widetilde{Y}^m_t,  \nonumber
\\
Z_t &:= \mathbb{I}_{\{t\leq \tau_m\}} \widetilde{Z}_t^1  + \sum_{m\geq 2} \mathbb{I}_{]\tau_{m-1}, \tau_m]}\widetilde{Z}^m_t, \nonumber
\\
N_t &:= \mathbb{I}_{\{t\leq \tau_m\}} \widetilde{N}_t^1  + \sum_{m\geq 2} \mathbb{I}_{]\tau_{m-1}, \tau_m]}\widetilde{N}^m_t.
\end{align*}
By (\ref{ae}), we have
$Y_{\cdot \wedge \tau_m} =\widetilde{Y}^m_{\cdot \wedge \tau_m}$, $\mathbb{I}_{\{t\leq \tau_m\}}Z_t = \mathbb{I}_{\{t\leq \tau_m\}}\widetilde{Z}^m_t$ and 
$N_{\cdot \wedge \tau_m} =\widetilde{N}^m_{\cdot \wedge \tau_m}$. Hence we can rewrite (\ref{bsdem}) as 
\[
Y_{t\wedge \tau_m} = \sup_n \inf_k Y^{n, k}_{\tau_m} +\int_{t\wedge \tau_m}^{\tau_m} \big( f{(s, Y_s, Z_s )}dA_s + g_s d\langle N\rangle_s \big) -\int_{t\wedge \tau_m}^{\tau_m}\big(Z_s dM_s + dN_s\big).
\]
 By sending $m$ to $ +\infty$, we prove that $(Y, Z, N)$ solves $(f, g, \xi)$. By (\ref{x1}),  we have 
\begin{align*}
|Y_t | =  |\sup_n \inf_k Y_t^{n, k} | \leq \frac{1}{\gamma} \ln \mathbb{E}\Big[ \exp\Big(\gamma e^{\beta A_{t, T}}|\xi| + 
\gamma \int_t^T e^{\beta A_{t, s}}{\alpha_s dA_s} \Big)\Big| \mathcal{F}_t \Big].
\end{align*}
\qed
\end{Proof}

Compared to    Mocha  and  Westray \cite{MW2012}, we prove the existence result under rather milder structure conditions. For example, \ref{as2'}(ii) gets rid of linear growth in $y$ and allows $g$ to be   any bounded process, which has been seen repeatedly throughout this paper. 
Secondly,  in contrast to their work, the assumption that $dA_t \ll c_A dt$, where $c_A$ is a positive constant, is not needed. 
Finally, 
they use a regularization procedure through quadratic BSDEs with bounded data. Hence, more demanding structure conditions are imposed to ensure that   the comparison theorem holds. On the contrary, the Lipschitz-quadratic regularization is more direct and essentially merely relies on \ref{as2'}
which is the most general assumption to our knowledge.
For the differences, the interested reader shall refer to \cite{M2009}, \cite{MW2012}.

Due to the same reason as in Proposition \ref{existence3}, the existence of a maximal or minimal solution is not available.

\begin{remm} Analogously to Hu and Schweizer \cite{HS2011}, 
one may easily extend the existence  result to infinite-horizon case. In abstract terms,  given exponential moments integrability on $\exp (\beta A_\infty)|\alpha|_\infty$, 
we regularize through Lipschitz-quadratic BSDEs with increasing horizons and null terminal value. Using a localization procedure  and the monotone stability result as in Theorem \ref{existence4}, we obtain a solution  which solves the infinite-horizon BSDE.
\end{remm}

As a result from Lemma \ref{aprioriestimate2}, we derive the estimates for the local martingale part. To save pages we only consider the following extremal case.
\begin{Corollary}[Estimate] \label{ubdestimate} Let   {\rm\ref{as2'}}  hold for $(f, g, \xi)$ and  $e^{\beta A_T} \big(|\xi|+ |\alpha|_T\big)$ has exponential moments of all orders. Then
any solution $(Y, Z, N)$ verifying {\rm(\ref{ey})} satisfies:
 $Y$ has exponential moments of all order and $Z\cdot M + N\in \mathcal{M}^p$ for all $p\geq 1.$ More precisely,  for all $p >1$, 
\begin{align*}
\mathbb{E}\big[e^{p\gamma Y^*} \big] \leq  \Big(\frac{p}{p-1}\Big)^p\mathbb{E}\Big[ \exp\Big( p\gamma e^{\beta A_T}\big(|\xi|+ |\alpha|_T \big)  \Big)\Big],
\end{align*}
and for all $p\geq 1$, 
\[
\mathbb{E}\Big[ \Big(\int_0^T\Big( Z_s^\top d\langle M\rangle_s Z_s + d\langle N\rangle_s \Big)\Big)^{\frac{p}{2}}    \Big]
\leq c \mathbb{E}\Big[ \exp\Big( 4p\gamma e^{\beta A_T}\big(|\xi| + |\alpha|_T\big)   \Big)    \Big], \]
where $c$ is a constant  only depending on $p, \gamma $.
\end{Corollary}
\begin{Proof}
The proof is exactly the same as Corollary 4.2,  Mocha and Westray \cite{MW2012} and hence omitted. \qed
\end{Proof}

Let us turn to the uniqueness result. 
We modify  Mocha and  Westray \cite{MW2012} to allow $g$ to be  any bounded  process rather than merely a constant. A convexity assumption  is imposed so as to use $\theta$-technique which proves to be convenient to treat quadratic terms. 
We start from comparison theorem and then move to uniqueness and stability result.  Similar results can be found in  Briand and  Hu \cite{BH2008} for Brownian setting
or  Da Lio and  Ley \cite{LL2006} from the point of view of PDEs.
To this end, the following structure conditions on $(f, g, \xi)$ are needed. 
\begin{as}\label{as3} There exist $\beta \geq 0, \gamma >0$  and an $\mathbb{R}^+$-valued $\Prog$-measurable process $\alpha$ such that $\mathbb{P}$-a.s.
\begin{enumerate}
\item [(i)] for any $t\in [0, T]$, $(y, z)\longmapsto f(t, y, z)$ is continuous;
\item [(ii)] $f$ is Lipschitz-continuous in $y$, i.e., for any $(t, z)\in [0, T]\times\mathbb{R}^d$, $y, y^\prime \in \mathbb{R}$, 
\[
|f(t, y, z) -f(t, y^\prime, z)| \leq \beta |y-y^\prime|;
\]
\item [(iii)] for any $(t, y)\in[0, T]\times\mathbb{R}$, $z\longmapsto f(t, y, z)$ is convex;
\item [(iv)] for any $(t, y, z)\in [0, T]\times\mathbb{R}\times\mathbb{R}^d$, 
\[|f(t, y, z)| \leq \alpha_t + \beta |y| + \frac{\gamma}{2}|\lambda_tz|^2.\]
\end{enumerate}
\end{as}

We start our proof of comparison theorem by observing that \ref{as3} implies \ref{as2'}. Hence existence  is ensured given suitable integrability.
Likewise, we keep in mind that $\mathbb{P}$-a.s. $ |g_\cdot| \leq \frac{\gamma}{2}$. 

\begin{thm}[Comparison Theorem]  
 \label{compare2}
Let $(Y, Z\cdot M + N)$, $(Y^\prime, Z^\prime\cdot  M + N^\prime)\in \mathcal{S}\times \mathcal{M}$ be solutions of
 $(f, g, \xi), (f^\prime, g^\prime, \xi^\prime)$, respectively, and $Y^*, (Y^\prime)^*$, $|\alpha|_T$ have exponential moments of all orders. If $\mathbb{P}$-a.s. for any $(t, y, z)\in [0, T]\times\mathbb{R}\times\mathbb{R}^d$, 
  $f(t, y, z)\leq f^\prime(t, y, z)$, $g_t \leq g^\prime_t$, $g_t^\prime \geq 0$, $\xi \leq \xi^\prime$ and $(f, g, \xi)$ verifies {\rm\ref{as3}}, 
then $\mathbb{P}$-a.s. $Y_\cdot\leq Y_\cdot^\prime$.
\end{thm}
\begin{Proof}
We introduce the notations used throughout the proof. 
For any $\theta \in (0, 1)$, define 
\begin{align*}
\delta f_t &:= f(t, Y_t^\prime, Z_t^\prime)-f^\prime(t, Y_t^\prime, Z_t^\prime),\\
 \delta_\theta Y &:= Y- \theta Y^\prime, \\
  \delta Y &:= Y-  Y^\prime,
\end{align*}
and  $\delta_\theta Z$,  $\delta Z$,  $\delta_\theta N, \delta N$, etc. analogously. Moreover, 
define \[
\rho_t : = \mathbb{I}_{\{ \delta_\theta Y_t \neq 0\}} \frac{f(t, Y_t, Z_t)-f(t, \theta Y^\prime_t, Z_t)}{\delta_\theta Y_t} .
\]
By  \ref{as3}(ii), $\rho$ is bounded by $\beta$ for any $\theta \in (0, 1)$. Hence $|\rho |_T \leq \beta \norm{A}$.
By It\^{o}'s formula, 
\begin{align*}
e^{|\rho|_t}\delta_\theta Y_t &= e^{|\rho|_T}\delta_\theta Y_T +\int_t^T e^{|\rho|_s}F^\theta_sdA_s +\int_t^T e^{|\rho|_s} \big(g_s d\langle N \rangle_s -\theta g^\prime_s d\langle N^\prime \rangle_s\big) \nonumber\\
&-\int_t^T e^{|\rho|_s}\big(\delta_\theta Z_s dM_s+ d\delta_\theta N_s\big), 
\end{align*}
where
\begin{align}
F^\theta_s &= f(s, Y_s, Z_s)- \theta f^\prime (s, Y_s^\prime, Z^\prime_s) -\rho_s \delta_\theta Y_s, \nonumber\\
  &=\theta \delta f_s + \big(f(s, Y_s, Z_s) - f(s, Y_s^\prime, Z_s)\big) + \big(f(s, Y_s^\prime, Z_s)- \theta f(s, Y^\prime_s, Z^\prime_s)\big)  -\rho_s \delta_\theta Y_s. \label{ftheta}
\end{align}
We then use \ref{as3}(ii)(iii) to deduce that 
\begin{align*}
f(s, Y_s, Z_s) - f(t, Y_s^\prime, Z_s)  &= f(s, Y_s, Z_s) - f(s, \theta Y_s^\prime, Z_s)  + f(s, \theta Y_s^\prime, Z_s) - f(s, Y_s^\prime, Z_s) \\
&= \rho_s \delta_\theta Y_s + f(t, \theta Y_s^\prime, Z_s) - f(s, Y_s^\prime, Z_s)  \\
&\leq \rho_s \delta_\theta Y_s + (1-\theta)\beta |Y_s^\prime|,\\
f(s, Y_s^\prime, Z_s)- \theta f(s, Y^\prime_s, Z^\prime_s) &= f(s, Y_s^\prime, \theta Z_t^\prime + (1-\theta)\frac{\delta_\theta Z_s}{1-\theta}) -\theta f (t, Y^\prime_s, Z^\prime_s) \\
& \leq (1-\theta) f(s, Y_s^\prime, \frac{\delta_\theta Z_s}{1-\theta})\\
&\leq (1-\theta)\alpha_s + (1-\theta)\beta |Y^\prime_s| + \frac{\gamma}{2(1-\theta)}|\lambda_s \delta_\theta Z_s|^2.
\end{align*}
We also note that $\mathbb{P}$-a.s. $\delta f_s \leq 0$.
Hence plugging these inequalities into (\ref{ftheta}) gives 
\begin{align}
F^\theta_s\leq
(1-\theta)\big( \alpha_s + 2\beta |Y_s^\prime| \big) + \frac{\gamma}{2(1-\theta)}|\lambda_s \delta_\theta Z_s|^2. \label{theta1}
\end{align}
We then perform an exponential transform to eliminate both quadratic terms. Set
\begin{align*}
c&:=  \frac{\gamma e^{\beta \norm{A}}}{1-\theta},\\
 P_t&: = \exp \big(ce^{|\rho|_t} \delta_\theta Y_t\big).
\end{align*}
By It\^{o}'s formula, 
\begin{align*}
P_t =P_T &+ \int_t^T  cP_s e^{|\rho|_s}\Big(
F_s^\theta -\frac{ce^{|\rho|_s}}{2}|\delta_\theta Z_s|^2
\Big)dA_s \\
& + \int_t^T cP_s e^{|\rho|_s}
 \Big(
g_s d\langle N \rangle_s -\theta g^\prime_s d\langle N^\prime \rangle_s - \frac{ce^{|\rho|_s}}{2}d\langle \delta_\theta N \rangle_s \Big)\\
&-\int_t^T cP_s e^{|\rho|_s} \big(\delta_\theta Z_s dM_s + d \delta_\theta N_s\big).
\end{align*}
For notational convenience, we define
\begin{align*}
G_t &:=cP_t e^{|\rho|_t}\Big(
F_t^\theta -\frac{ce^{|\rho|_t}}{2}|Z_t^\theta|^2
\Big),\\
H_t &:= \int_0^t cP_s e^{|\rho|_s}
 \Big(
g_s d\langle N \rangle_s -\theta g^\prime_s d\langle N^\prime \rangle_s - \frac{ce^{|\rho|_s}}{2}d\langle N^\theta \rangle_s \Big).
\end{align*}
By (\ref{theta1}), we have
\[
G_t = cP_te^{|{\rho}|_t} \Big( (1-\theta)\big(\alpha_t + 2\beta |Y^\prime_t| \big) \Big) \leq P_tJ_t, 
\]
where \[
J_t: = \gamma e^{2\beta \norm{A}}\big( \alpha_t + 2\beta |Y_t^\prime|\big).
\]
We claim that   $H$ can also be eliminated. Indeed, 
\begin{align*}
d \langle \delta_\theta N \rangle & = d\langle N \rangle + \theta^2 d\langle N^\prime \rangle -2\theta d \langle N, N^\prime \rangle\\
&\gg d\langle N\rangle + \theta^2 d\langle N^\prime \rangle -\theta  d\langle N\rangle  -\theta d\langle N^\prime\rangle \\
& =(1-\theta)\big(d\langle N \rangle - \theta d\langle N^\prime \rangle\big)\\
& =  (1-\theta)  d \delta_\theta \langle N \rangle.
\end{align*}
We then come back to $H$ and use this inequality to deduce that   
\begin{align*}
g_t d\langle N \rangle_t -\theta g^\prime_t d\langle N^\prime \rangle_t - \frac{ce^{|{\rho}|_t}}{2}d\langle \delta_\theta N \rangle_t &= 
g^+_td\langle N\rangle_t - g^-_t d\langle N\rangle_t -\theta  g_t^\prime  d\langle N^\prime \rangle_t -\frac{ce^{|\rho|_t}}{2}d\langle \delta_\theta N \rangle_t \\
&\ll g^+_t d \delta_\theta \langle N\rangle_t + \theta (g^+_t - g^\prime_t) d\langle N^\prime \rangle_t  -\frac{ce^{|{\rho}|_t}}{2}d\langle \delta_\theta N\rangle_t \\
&\ll g^+_td\delta_\theta\langle N\rangle_t -\frac{\gamma}{2(1-\theta)}d\langle \delta_\theta N \rangle_t \\
&\ll 0,
\end{align*}
 due to $g_\cdot^+ \leq g_\cdot^\prime$ and $g_\cdot\leq \frac{\gamma}{2}$. Hence $dH_t \ll 0.$
To eliminate $G$, we set
$
D_t: = \exp\big( |J|_t \big).
$
By It\^o's formula, 
\begin{align*}
d(D_tP_t) &= D_t \Big( \big(P_tJ_t -G_t \big)dA_t -dH_t + cP_t e^{|{\rho}|_t }\big(\delta_\theta Z_tdM_t +d\delta_\theta N_t \big)  \Big).
\end{align*}
But  previous results show that 
$(P_tJ_t -G_t)dA_t - dH_t \gg 0$. Hence $DP$ is a local submartingale.
 Thanks to the exponential moments integrability on $|\alpha|_T$ and $(Y^\prime)^*$ (and hence $|J|_T$), 
we use a localization procedure and easily deduce that 
\begin{align}
P_t &\leq \mathbb{E}\Big[\exp\Big( \int_t^{T} J_s dA_s \Big)P_{T}\Big| \mathcal{F}_t  \Big]. \label{psub}
\end{align}
We  come back to the definition of $P_T$ and observe that 
\begin{align*}
\delta_\theta \xi &\leq (1-\theta) |\xi|+ \theta\delta \xi\\
 &\leq (1-\theta) |\xi|.
\end{align*}
Hence (\ref{psub}) gives
\begin{align*}
\exp\Big( \frac{\gamma e^{\beta \norm{A} + |{\rho}|_t}}{1-\theta}\delta_\theta Y_t \Big)& \leq 
\mathbb{E}\Big[ \exp\Big( \int_t^T J_s      dA_s      \Big) 
\exp\big( ce^{|{\rho}|_T}\delta_\theta \xi	      \big) \Big| \mathcal{F}_t
\Big]  \\ 
&\leq
\mathbb{E}\Big[ \exp\Big( \int_t^T J_s   dA_s      \Big) 
\exp\big( \gamma e^{2\beta \norm{A}}|\xi|	      \big) \Big| \mathcal{F}_t
\Big].
\end{align*}
Hence
\begin{align*}
\delta_\theta Y_t \leq \frac{1-\theta}{\gamma}\ln \mathbb{E}\Big[\exp \Big(\gamma e^{2\beta \norm{A} }\Big( |\xi| + \int_t^T \big(\alpha_s + 2\beta |Y_s^\prime|\big) dA_s    \Big)            \Big) \Big| \mathcal{F}_t  \Big].
\end{align*}
Therefore we obtain $\mathbb{P}$-a.s.  $Y_t \leq Y_t^\prime$, by sending $\theta$ to $1$. By  the continuity of $Y$ and $Y^\prime$, we also have $\mathbb{P}$-a.s. $Y_\cdot\leq Y_\cdot^\prime$.
\qed
\end{Proof}
As a byproduct, we can prove the existence of a unique solution given \ref{as3}. 
\begin{Corollary}[Uniqueness]
\label{exunique}
If $(f, g, \xi)$ satisfies {\rm\ref{as3}}, $\mathbb{P}$-a.s. $g_\cdot\geq 0$ and $|\xi|$, $|\alpha|_T$ have exponential moments of all orders, then there exists a unique solution $(Y, Z, N)$ to $(f, g, \xi)$ such that $Y^*$ has exponential moments of all order and $(Z\cdot M + N)\in \mathcal{M}^p$ for all $p\geq 1$.
\end{Corollary}
\begin{Proof}
The existence  of a unique solution in the above sense is immediate from  Theorem \ref{existence4} (existence), Theorem \ref{compare2} (comparison theorem) and Corollary \ref{ubdestimate} (estimate).
\qed
\end{Proof}
\begin{remm}
There are spaces  to sharpen the uniqueness.   The convexity in $z$ motivates one to replace \ref{as3}(iv) by 
\[
-\underline{\alpha}_t - \beta |y| - \kappa |\lambda_t z|
\leq f (t, y, z) \leq \overline{\alpha}_t +\beta |y| + \frac{\gamma}{2}|\lambda_t z|^2.
\]
Secondly, in  view of Delbaen et al \cite{DHR2011}, 
  we may prove uniqueness given weaker integrability, by characterizing
the  solution as the value process of a stochastic control problem. 
\end{remm}

It turns out that a stability result also holds given convexity condition. The proof is a modification of Theorem \ref{compare2} (comparison theorem). We set $\mathbb{N}^0:= \mathbb{N}^+ \cup \{ 0 \}.$
\begin{prop}[Stability] 
Let $(f^n, g^n, \xi^n)_{n\in \mathbb{N}^0}$  with $g^n_\cdot \geq 0$ $\mathbb{P}$-a.s. 
satisfy {\rm\ref{as3}} associated with $(\alpha^n, \beta, \gamma, \varphi)$, and $(Y^n, Z^n, N^n)$ be their unique solutions in the sense of Corollary {\rm\ref{exunique}}, respectively.
If $\xi^n - \xi^0 \longrightarrow 0$, 
$\int_0^T |f^n -f^0 | (s, Y_s^0, Z_s^0) dA_s  \longrightarrow 0$ in probability, $\mathbb{P}$-a.s. $g^n_\cdot -g^0_\cdot \longrightarrow 0$  as $n$ goes to $+\infty$ and for each $p > 0$, 
\begin{align}
&\sup_{n\in \mathbb{N}^0} \mathbb{E}\Big[ \exp\Big(   p \big(|\xi^n\big| + |\alpha^n|_T \big)           \Big)\Big] <+\infty,   \label{bdddata}
 \\
& \sup_{n\in\mathbb{N}^0} |g^n_\cdot| \leq  \frac{\gamma}{2}\ \mathbb{P}\text{-a.s.}  \nonumber
\end{align}
Then  for each $p \geq  1$, 
\begin{align*}
&\lim_n \mathbb{E}\Big[ \exp\big( p|Y^n-Y^0|^* \big)\Big]  =1, \\
&\lim_n \mathbb{E}\Big[  \Big(\int_0^T \Big((Z_s^n -Z_s^0)^\top d\langle M\rangle_s (Z_s^n -Z_s^0)  +d\langle N^n- N^0 \rangle_s   \Big) \Big)^{\frac{p}{2}} \Big] = 0.
\end{align*}
\end{prop}
\begin{Proof} By Corollary \ref{ubdestimate} (estimate), for any $p \geq 1$, 
\begin{align}
\sup_{n\in \mathbb{N}^0}\mathbb{E}\Big[ \exp\big(p(Y^n)^*\big) + \Big(\int_0^T \Big( (Z_s^n)^\top d\langle M \rangle_s Z_s^n + d\langle N^n\rangle_s \Big) \Big)^{\frac{p}{2}}        \Big] <+\infty. \label{bddn}
\end{align}
Hence the sequence of random variables
\begin{align*}
\exp \Big(p|Y^n- Y^0|^* \Big)  + \Big(\int_0^T \Big((Z_s^n - Z^0_s)^\top d\langle M \rangle_s (Z_s^n - Z^0_s) + d\langle N^n -N^0 \rangle_s\Big)\Big)^{\frac{p}{2}}
\end{align*}
is uniformly integrable. Due to Vitali  convergence,  it is hence sufficient to prove that 
\begin{align*}
|Y^n - Y|^* + \int_0^T \Big((Z_s^n -Z_s^0)^\top  d\langle M\rangle (Z_s^n -Z_s^0) + d\langle N^n - N\rangle_s\Big) \longrightarrow 0
\end{align*}
in probability as $n$ goes to $+\infty$. 

(i).  We prove $u.c.p$ convergence of $Y^n - Y^0$. To this end we use $\theta$-technique in  the spirit of Theorem \ref{compare2} (comparison theorem). For any $\theta \in (0, 1)$, define 
\begin{align*}
\delta f_t^n &: = f^0(t, Y_t^0, Z_t^0)  - f^n (t, Y_t^0, Z_t^0),\\
\delta g^n &: = g^0- g^n,\\
 \delta_\theta Y^n &: = Y^0- \theta Y^n, 
\end{align*}
 and  $\delta_\theta Z^n, \delta_\theta N^n$, $\delta_\theta \langle N\rangle^n$, etc. analogously.  Further, set
\begin{align*}
\rho_t &: =\mathbb{I}_{\{ Y_t^0 - Y_t^n \neq 0\}} \frac{f^n (t, Y_t^0, Z^n_t) -f^n (t, Y_t^n, Z_t^n)}{Y_t^0 - Y_t^n}, \\
c&: = \frac{\gamma e^{\beta \norm{A}}}{1-\theta}, \\
P_t^n &: = \exp \big(c e^{|\rho|_t} \delta_\theta Y_t^n \big), \\
 J_t^n&: = \gamma e^{2  \beta \norm{A}}\big(\alpha_t^n + 2\beta |Y^0_t|\big), \\
 D_t^n &: = \exp \Big(\int_0^t J_s^n dA_s \Big).
\end{align*} 
Obviously $\rho$ is bounded by $\beta $ due to \ref{as3}(i).
 The $\theta$-difference implies that 
\begin{align}
&f^0(t, Y^0_t, Z^0_t) - \theta f^n (t, Y_t^n, Z_t^n) \nonumber \\
& = \delta f_t^n  + \big( \theta f^n(t, Y_t^0, Z_t^n)  - \theta f^n(t, Y^n_t, Z^n_t) \big) +  \big(f^n (t, Y^0_t, Z^0_t) - \theta f^n(t, Y_t^0, Z_t^n)\big). \label{0n}
\end{align}
By \ref{as3}(i)(ii), 
\begin{align*}
\theta f^n(t, Y_t^0, Z_t^n) - \theta f^n(t, Y^n_t, Z^n_t)
& =\theta \rho_t (Y_t^0 - Y_t^n)\\
& = \rho_t \big(\theta Y_t^0  - Y_t^0 + Y_t^0	 -\theta Y_t^n\big)\\
&\leq 
 (1-\theta)\beta |Y_t^0| + \rho_t \delta_\theta Y^n_t,\\
f^n (t, Y^0_t, Z^0_t) - \theta f^n(t, Y_t^0, Z_t^n) &\leq  (1-\theta)\alpha_t^n + (1-\theta) \beta |Y_t^0| + \frac{\gamma}{2(1-\theta)} |\delta_\theta Z^n_t|^2.
\end{align*}
Hence  (\ref{0n}) gives
\begin{align}
f^0(t, Y^0_t, Z^0_t) - \theta f^n (t, Y_t^n, Z_t^n)  -\rho_t \delta_\theta Y^n_t \leq \delta f_t^n + (1-\theta) \big(\alpha_t^n + 2 \beta |Y_t^0| \big) +\frac{\gamma}{2(1-\theta)}|\delta_\theta Z_t^n|^2. \label{thetaf}
\end{align}
To analyze the quadratic term concerning $N^0$ and $N^n$, we deduce by the same arguments as in Theorem \ref{compare2}  that 
\begin{align}
g^0_t d\langle N^0 \rangle_t -\theta g^n_t d\langle N^n \rangle_t - \frac{ce^{|\rho|_t}}{2}d\langle 
\delta_\theta N \rangle_t
&= \delta g_t^n d\langle N^0\rangle_t + g^n_t d\delta_\theta\langle N \rangle_t^n  -\frac{ce^{|\rho|_t}}{2}d\langle \delta_\theta N^n \rangle_t \nonumber\\
&\ll g^n_t \Big(d\delta_\theta\langle N\rangle^n_t -\frac{1}{1-\theta}d\langle \delta_\theta N^n\rangle_t \Big) +  \delta g^n_t d\langle N^0 \rangle_t \nonumber\\
&\ll \delta g^n_t d\langle N^0 \rangle_t. \label{thetag}
\end{align}
Given (\ref{thetaf}) and (\ref{thetag}), we use an  exponential transform which is analogous to that in Theorem \ref{compare2}. This  gives
\begin{align*}
P_t^n \leq D_t^n P_t^n \leq \mathbb{E}\Big[ D_T^n P_T^n + \frac{\gamma e^{2\beta \norm{A}}}{1-\theta}\int_t^T D_s^n P_s^n\big(| \delta f_s^n| dA_s       
+ | \delta g_s^n |d\langle N^0\rangle_s  \big)
\Big| \mathcal{F}_t\Big].
\end{align*}
Using $\log x \leq x$ and $Y^0 - Y^n \leq (1-\theta)|Y^n| + \delta_\theta Y^n$, we deduce that 
\[
Y_t^0 - Y_t^n \leq (1-\theta) |Y_t^n| + \frac{1-\theta}{\gamma} \mathbb{E}\Big[ D_T^n P_T^n + \frac{\gamma e^{2\beta \norm{A}}}{1-\theta}\int_t^T D_s^n P_s^n \big(|\delta f_s^n| dA_s       
+  |\delta g_s^n| d\langle N^0\rangle_s\big)  
\Big| \mathcal{F}_t\Big].
\]
Set
\begin{align*}
\Lambda^n(\theta) &: = \exp\Big( \frac{\gamma e^{2\beta \norm{A}}}{1-\theta} \big((Y^0)^* + (Y^n)^*\big)     \Big) \geq P_t^n, \\
\Xi^n (\theta) &:=  \exp\Big( \frac{\gamma e^{2\beta \norm{A}}}{1-\theta} \big(|\xi^0 -\theta \xi^n| \vee |\xi^n-\theta \xi^0|\big)     \Big) \geq P_T^n.
\end{align*}
We then have
\[
Y_t^0 - Y_t^n \leq (1-\theta) |Y_t^n| + \frac{1-\theta}{\gamma} \mathbb{E}\Big[ D_T^n\Xi^n(\theta)  + \frac{\gamma e^{2\beta \norm{A}}}{1-\theta}D_T^n\Lambda^n(\theta)\int_t^T \big(| \delta f_s^n| dA_s       
+   |\delta g_s^n| d\langle N^0\rangle_s  \big)
\Big| \mathcal{F}_t\Big].
\]
Now we use \ref{as3}(ii)(iii) to $f^n$ and proceed analogously to Theorem \ref{compare2}. This gives
\[
Y_t^n -Y_t^0 \leq (1-\theta)|Y_t^0| +   \frac{1-\theta}{\gamma} \mathbb{E}\Big[ D_T^n \Xi^n(\theta) + \frac{\gamma e^{2\beta  \norm{A}}}{1-\theta}D_T^n\Lambda^n(\theta)\int_t^T \big(| \delta f_s^n| dA_s       
+    |\delta g_s^n| d\langle N^0\rangle_s \big)  
\Big| \mathcal{F}_t\Big].
\]
Though looking symmetric, the two inequalities come from slightly different treatments for the $\theta$-difference. The two estimates give 
\begin{align*}
|Y_t^n -Y_t^0| \leq \underbrace{(1-\theta)\big(|Y_t^0| +|Y_t^n|\big)}_\text{$X^1_t$}  &+ \underbrace{\frac{1-\theta}{\gamma}\mathbb{E}\Big[ D_T^n \Xi^n(\theta)\Big|\mathcal{F}_t\Big]}_{\text{$X^2_t$}} \\
&+ 
\underbrace{e^{2\beta \norm{A}}\mathbb{E}\Big[D_T^n\Lambda^n(\theta)\int_0^T \big(| \delta f_s^n| dA_s     +  |\delta g_s^n| d\langle N^0\rangle_s \big)  
\Big| \mathcal{F}_t\Big]}_{\text{$X^3_t$}}.
\end{align*}
We then prove $u.c.p$ convergence of $Y^n-Y^0.$ For any $\epsilon >0$, 
\begin{align}
\mathbb{P}\Big( |Y^n -Y^0|^* \geq  \epsilon \Big)\leq \mathbb{P}\Big( (X^1)^* \geq \frac{\epsilon}{3}    \Big) 
+\mathbb{P}\Big( (X^2)^* \geq \frac{\epsilon}{3}    \Big) 
+\mathbb{P}\Big( (X^3)^* \geq \frac{\epsilon}{3}    \Big). \label{converge}
\end{align}
We aim at showing that each term on the right-hand side of (\ref{converge}) converges to $0$ if we send $n$ to $+\infty$ first and then $\theta$ to $1$.
To this end, we give some useful estimates.
By Chebyshev's inequality,  
\begin{align*}
\mathbb{P}\Big( (X^1)^* \geq \frac{\epsilon}{3}    \Big) \leq \frac{3(1-\theta)}{\epsilon}\mathbb{E}\big[ (Y^0)^* + (Y^n)^* \big],
\end{align*}
where  $\mathbb{E}[ (Y^0)^* + (Y^n)^* ]$ is uniformly bounded.
Secondly, 
Doob's inequality yields
\begin{align}
\mathbb{P}\Big( (X^2)^* \geq \frac{\epsilon}{3}    \Big)  \leq  \frac{3(1-\theta)\gamma}{\epsilon}\mathbb{E}\big[D_T^n    \Xi_T^n    \big]. \label{x3}
\end{align}
Moreover, by Vitali convergence,  the right-hand side of (\ref{x3}) satisfies
\begin{align*}
\limsup_n \mathbb{E}\big[D_T^n \Xi_T^n \big] &\leq \sup_n\mathbb{E}\big[(D^n)^2\big]^{\frac{1}{2}}\cdot 
\limsup_n\mathbb{E}\big[ (\Xi^n)^2\big]^{\frac{1}{2}}
\\
&\leq \sup_n\mathbb{E}\big[(D^n)^2\big]^{\frac{1}{2}}\cdot 
\mathbb{E}\Big[ \exp\Big(2 \gamma e^{2\beta \norm{A}}|\xi^0|\Big)\Big]^{\frac{1}{2}}\\
 & < +\infty.
\end{align*}
Hence, the first term and the second term on the right-hand side of (\ref{converge}) converge to $0$ as $n$ goes to $+\infty$ and $\theta$  goes to $1$.
Finally, we claim that the third term on the right-hand side of (\ref{converge}) also converges. Indeed,  Doob's inequality and H\"{o}lder's inequality give
\begin{align}
\mathbb{P}\Big( (X^3)^* \geq \frac{\epsilon}{3}    \Big) &\leq \frac{3e^{2\beta\norm{A}}}{\epsilon}\mathbb{E}\Big[D_T^n\Lambda^n(\theta)\int_t^T \big(| \delta f_s^n| dA_s     +  |\delta g_s^n| d\langle N^0\rangle_s \big)  
\Big] \nonumber \\
&\leq \frac{3e^{2\beta\norm{A}}}{\epsilon}\mathbb{E}\Big[\big(D_T^n\Lambda^n(\theta) \big)^2 \Big]^{\frac{1}{2}}
\mathbb{E}\Big[\Big(\int_0^T \big(   | \delta f_s^n| dA_s     +  |\delta g_s^n| d\langle N^0\rangle_s \big) \Big)^{2}    \Big]^{\frac{1}{2}}. \label{x3vitali}
\end{align}
Note that 
\[
\int_0^T \big(   | \delta f_s^n| dA_s     +  |\delta g_s^n| d\langle N^0\rangle_s \big)
\leq |\alpha|_T + |\alpha^n|_T 
+ 2\norm{A} (Y^0)^*  + \gamma \langle Z^0\cdot M + N^0\rangle_T.
\]
Hence the left-hand side of this inequality has finite moments of all orders by Corollary \ref{ubdestimate}. 
Therefore, the left-hand side of (\ref{x3vitali}) converges to $0$ as $n$ goes to $+\infty$ due to Vitali convergence. 

Finally, collecting these convergence results for each term in (\ref{converge})  gives the convergence of $Y^n-Y^0$.

(ii).  It remains to prove convergence of the martingale parts. 
By It\^{o}'s formula, 
\begin{align*}
&\mathbb{E}\Big[ \int_0^T \Big( (Z_s^n - Z_s^0)^\top d\langle M \rangle_s (Z_s^n - Z_s^0) + d\langle N^n -N^0\rangle_s \Big)  \Big]\\
&\leq
  \mathbb{E}\big[ \big|\xi^n - \xi^0 \big|^2\big]  + 2\mathbb{E}\Big[ |Y^n-Y^0|^* \int_0^T \big|F^n(s, Y_s^n, Z_s^n) - F^0(s, Y^0_s, Z^0_s)\big| dA_s \Big]
\\
&+ 2\mathbb{E}\Big[ |Y^n-Y^0|^* \Big| \int_0^T \big(g_s^nd\langle N^n \rangle_s - g^0_s d\langle N^0 \rangle_s \big)         \Big|    \Big],
\end{align*}
As before, we conclude by Vitali convergence. 
\qed
\end{Proof}

\subsection{Change of Measure}
\label{section24}
In the final section, we show that given exponential moments integrability, the martingale part $Z\cdot M + N$, though not BMO,   defines an equivalent change of measure, i.e., its stochastic exponential is a  strictly positive martingale. 
We don't require convexity which ensures uniqueness. But to derive the estimate for
$\int_0^T f(s, Y_s, Z_s)dA_s$, we use \ref{as2''} where
$f$ is of linear growth in $y$. 
We keep assuming that $\mathbb{P}$-a.s. $|g_\cdot|\leq \frac{\gamma}{2}$. The following result comes from Mocha and Westray \cite{MW2012}. 
\begin{thm}[Change of Measure]
If $(f, g, \xi)$ satisfies {\rm \ref{as2''}} and $\xi$, $|\alpha|_T$ have exponential moments of all orders, then for any solution $(Y, Z, N)$ such that $Y$ has exponential moments of all orders and any
 $|q| > \frac{\gamma}{2}$, $\mathcal{E}\big(q\big(Z\cdot M + N\big)\big)$ is a continuous martingale.
\end{thm}
\begin{Proof}
We start by recalling  Lemma 1.6. and Lemma 1.7., Kazamaki \cite{K1994}: if $\widetilde{M}$ is a martingale such that 
\begin{align}
\sup_{\tau \in \mathcal{T}}\mathbb{E}\Big[
\exp\Big(\eta \widetilde{M}_\tau      + \Big(\frac{1}{2}- \eta\Big)\langle \widetilde{M} \rangle_\tau \Big)
\Big]  < +\infty, \label{criteria}
\end{align} 
for  $\eta \neq 1$, then $\mathcal{E}\big(\eta \widetilde{M}\big)$ is a martingale. Moreover, if (\ref{criteria}) holds for some $\eta^* > 1$ then it holds for any $\eta \in (1, \eta^*)$.

By  Lemma \ref{ubdestimate} (estimate),  $Z\cdot M + N$ is a continuous martingale. 
First of all, we apply the above criterion to $\widetilde{M}: = \tilde{q}(Z\cdot M + N)$
 for some fixed $|\tilde{q}| > \frac{\gamma}{2}$. Define $ \Lambda_t (\eta)$ such that
 \[
 \ln \Lambda_t (\eta) : = 
 \tilde{q}\eta \big( 
(Z\cdot M)_t +N_t\big)  + \tilde{q}^2 \Big(\frac{1}{2} - \eta\Big)\langle Z\cdot M + N \rangle_t. 
 \]
From  the BSDE (\ref{sbsde})
 and \ref{as2''}, we obtain, for any $\tau \in \mathcal{T}$,
\begin{align}
\ln \Lambda_\tau (\eta)
&=
\tilde{q}\eta \Big(
Y_t - Y_0 +\int_0^t \big(f(s, Y_s, Z_s)dA_s + g_s d\langle N\rangle_s
\big) \Big)
+\tilde{q}^2 \Big( \frac{1}{2} -\eta\Big)\langle Z\cdot M  + N\rangle_t \nonumber\\
&\leq
(2 + \beta\norm{A})|\tilde{q}|\eta  Y^* + |\tilde q|\eta |\alpha|_T  + 
|\tilde{q}|\eta \Big( \frac{\gamma}{2} +\frac{|\tilde{q}|}{\eta}\Big( \frac{1}{2}-\eta\Big)   \Big) \langle Z\cdot M + N \rangle_T. \label{lnlambda}
\end{align}
Note that 
\[
 \frac{\gamma}{2} +\frac{|\tilde{q}|}{\eta}\Big( \frac{1}{2}-\eta\Big)  \leq 0  \Longleftrightarrow \eta \geq \frac{|\tilde{q}|}{2|\tilde{q}|-\gamma} =: q_0 \Big(>\frac{1}{2}\Big).
\]
Hence for any $\eta \geq q_0$, (\ref{lnlambda}) gives
\[
\Lambda_\tau (\eta) \leq 
\exp \big( |\tilde{q}|\eta (2 + \beta) Y_* + |\tilde q|\eta |\alpha|_T    \big).
\]
By exponential moments integrability, we have
\[
\sup_{\tau \in \mathcal{T}}
\mathbb{E}\big[ \Lambda_\tau (\eta)  \big] < + \infty.
\]
It then follows from the first statement of the criterion that $\mathcal{E}\big(\tilde{q}\eta (Z\cdot M + N)\big)$ is a martingale for all $\eta \in [q_0, \infty)\backslash\{1\}$. The second statement ensures that it is a martingale for any $\eta >1$. For any $|q| > \frac{\gamma}{2}$, we set 
 $|\tilde{q}|\in (\frac{\gamma}{2},|q|)$, $\eta := \frac{q}{\tilde{q}} > 1$, and apply the result above to conclude that  $\mathcal{E}\big(q \big(Z\cdot M+ N \big)\big)$ is a martingale.
\qed
\end{Proof}
\textbf{Acknowledgement.} The author thanks Martin Schweizer for his supervision and many helpful remarks.
\bibliography{myrefs}
\bibliographystyle{plain}
\end{document}